\documentstyle[txmac,a4,
amssymb,%
case,%
twoside,%
nocaphead,%
epsf,%
mypic,times,mathptm]{article}

\advance\oddsidemargin by -1.8cm
\advance\evensidemargin by -1.8cm
\advance\textwidth by 3.6cm

\def\mynewtheo#1#2{%
\newtheorem{@#1}{#2}
\newenvironment{#1}{\begin{@#1}\rm}{\end{@#1}}}

\mynewtheo{lemma}{Lemma}
\mynewtheo{theo}{Theorem}
\mynewtheo{rem}{Remark}
\mynewtheo{defi}{Definition}
\mynewtheo{conj}{Conjecture}
\mynewtheo{corr}{Corollary}
\mynewtheo{prop}{Proposition}
\mynewtheo{question}{Question}
\mynewtheo{exam}{Example}

\newenvironment{theorem}{\begin{theo}}{\end{theo}}

\begin{document}

\makeatletter

\newenvironment{eqn}{\begin{equation}}{\end{equation}\@ignoretrue}

\newenvironment{myeqn*}[1]{\begingroup\def\@eqnnum{\reset@font\rm#1}%
\xdef\@tempk{\arabic{equation}}\begin{equation}\edef\@currentlabel{#1}}
{\end{equation}\endgroup\setcounter{equation}{\@tempk}\ignorespaces}

\newenvironment{myeqn}[1]{\begingroup\let\eq@num\@eqnnum
\def\@eqnnum{\bgroup\let\r@fn\normalcolor 
\def\normalcolor####1(####2){\r@fn####1#1}%
\eq@num\egroup}%
\xdef\@tempk{\arabic{equation}}\begin{equation}\edef\@currentlabel{#1}}
{\end{equation}\endgroup\setcounter{equation}{\@tempk}\ignorespaces}

\newenvironment{myeqn**}{\begin{myeqn}{(\arabic{equation})\es\es\mbox{\qed}}\edef\@currentlabel{\arabic{equation}}}
{\end{myeqn}\stepcounter{equation}}

\def\bysame{\same[\kern2cm]\,}
\def\lfra{\leftrightarrow}
\def\qed{\hfill\@mt{\Box}}
\def\@mt#1{\ifmmode#1\else$#1$\fi}
\def\qqed{\hfill\@mt{\Box\enspace\Box}}

\let\ap\alpha
\let\tl\tilde
\let\sg\sigma
\let\dl\delta
\let\Dl\Delta
\let\eps\varepsilon
\let\gm\gamma
\let\bt\beta
\let\fa\forall
\let\nb\nabla

\def\cF{{\cal F}}
\def\bZ{{\Bbb Z}}
\def\bN{{\Bbb N}}
\def\spn{\mbox{\operator@font span}\,}
\def\int{\mbox{\operator@font int}\,}
\def\ext{\mbox{\operator@font ext}\,}
\def\val{\mbox{\operator@font val}\,}
\let\ds\displaystyle
\def\cf{\text{\rm cf}\,}
\def\mcf{\min\cf\,}
\def\mc{\min\cf\,}
\def\md{\min\deg}
\def\Md{\max\deg}
\let\es\enspace
\let\sm\setminus
\let\wh\widehat
\def\hG{\hat G}
\def\nin{\not\in}
\def\br#1{\left\langle#1\right\rangle}
\def\lra{\longrightarrow}
\def\Ra{\Rightarrow}

\def\ssim{\mathrel{\raise-0.3em\hbox{$\stackrel{\ds \sim}{\vbox{\vskip-0.2em\hbox{$\scriptstyle *$}}}$}}}

\def\ellipse#1#2#3{
  \pictranslate{#1}{
    \picrotate{#2}{
      \picellipse{0 0}{#3}{}
    }
  }
}

\def\drop#1#2#3{
  \pictranslate{#1}{
    \picfillgraycol{0}
    \picrotate{#2}{
      \picfilledcircle{0 0}{0.08}{}
      \picscale{#3}{
        \picmultigraphics[S]{2}{-1 1}{
          \piccurve{0 0}{0.5 d}{0.5 1}{0 1}
	}
      }
    }
  }
}  
\def\ddrop#1#2#3{{\piclinedash{0.15 0.05}{0.25}\drop{#1}{#2}{#3}}}

\def\vrt#1{{\picfillgraycol{0}\picfilledcircle{#1}{0.08}{}}}
\def\cycl#1#2#3#4{\vrt{#1}\vrt{#2}\vrt{#3}\vrt{#4}%
\picline{#1}{#2}\picline{#2}{#3}\picline{#3}{#4}\picline{#4}{#1}%
}

\def\@curvepath#1#2#3{%
  \@ifempty{#2}{\piccurveto{#1 polar}{@stc}{@std}#3}%
    {\piccurveto{#1 polar}{#2 polar}{#2 polar #3 polar 0.5 conv}
    \@curvepath{#3}}%
}
\def\curvepath#1#2#3{%
  \piccurve{#1 polar}{#2 polar}{#2 polar}{#2 polar #3 polar 0.5 conv}%
  \picPSgraphics{/@stc [ #1 polar #2 polar -1 conv ] $ D /@std [ #1 polar ] $ D }%
  \@curvepath{#3}%
}

\def\@opencurvepath#1#2#3{%
  \@ifempty{#3}{\piccurveto{#1}{#1}{#2}}%
    {\piccurveto{#1}{#2}{#2 #3 0.5 conv}\@opencurvepath{#3}}%
}
\def\opencurvepath#1#2#3{%
  \piccurve{#1}{#2}{#2}{#2 #3 0.5 conv}%
  \@opencurvepath{#3}%
}

\def\epsfs#1#2{{\catcode`\_=11\relax\ifautoepsf\unitxsize#1\relax\else
\epsfxsize#1\relax\fi\epsffile{#2.eps}}}
\def\epsfsv#1#2{{\vcbox{\epsfs{#1}{#2}}}}
\def\vcbox#1{\setbox\@tempboxa=\hbox{#1}\parbox{\wd\@tempboxa}{\box
  \@tempboxa}}
\def\lz{\linebreak[0]\verb}

\def\@test#1#2#3#4{%
  \let\@tempa\go@
  \@tempdima#1\relax\@tempdimb#3\@tempdima\relax\@tempdima#4\unitxsize\relax
  \ifdim \@tempdimb>\z@\relax
    \ifdim \@tempdimb<#2%
      \def\@tempa{\@test{#1}{#2}}%
    \fi
  \fi
  \@tempa
}

\def\go@#1\@end{}
\newdimen\unitxsize
\newif\ifautoepsf\autoepsftrue

\unitxsize4cm\relax
\def\epsfsize#1#2{\epsfxsize\relax\ifautoepsf
  {\@test{#1}{#2}{0.1 }{4   }
		{0.2 }{3   }
		{0.3 }{2   }
		{0.4 }{1.7 }
		{0.5 }{1.5 }
		{0.6 }{1.4 }
		{0.7 }{1.3 }
		{0.8 }{1.2 }
		{0.9 }{1.1 }
		{1.1 }{1.  }
		{1.2 }{0.9 }
		{1.4 }{0.8 }
		{1.6 }{0.75}
		{2.  }{0.7 }
		{2.25}{0.6 }
		{3   }{0.55}
		{5   }{0.5 }
		{10  }{0.33}
		{-1  }{0.25}\@end
		\ea}\ea\epsfxsize\the\@tempdima\relax
		\fi
		}

\let\old@tl\~
\def\~{\raisebox{-0.8ex}{\tt\old@tl{}}}

\author{A. Stoimenow\footnotemark[1]\\[2mm]
\small Department of Mathematics, \\
\small University of Toronto,\\
\small Canada M5S 3G3\\
\small e-mail: {\tt stoimeno@math.toronto.edu}\\
\small WWW: {\hbox{\tt http://www.math.toronto.edu/stoimeno/}}
}

{\def\thefootnote{\fnsymbol{footnote}}
\footnotetext[1]{Supported by a DFG postdoc grant.}
}

\title{\large\bf \uppercase{On polynomials and surfaces of
variously positive links}\\[4mm]
{\small\it This is a preprint. 
I would be grateful for any comments and corrections!}}

\date{\large Current version: \curv\ \ \ First version:
\makedate{4}{1}{2002}}

\def\Myfrac#1#2{\mbox{\small$\ds\frac{#1}{#2}$}}

\maketitle

\long\def\@makecaption#1#2{%
   \vskip 10pt
   {\let\label\@gobble
   \let\ignorespaces\@empty
   \xdef\@tempt{#2}%
   }%
   \ea\@ifempty\ea{\@tempt}{%
   \setbox\@tempboxa\hbox{%
      \fignr#1#2}%
      }{%
   \setbox\@tempboxa\hbox{%
      {\fignr#1:}\capt\ #2}%
      }%
   \ifdim \wd\@tempboxa >\captionwidth {%
      \rightskip=\@captionmargin\leftskip=\@captionmargin
      \unhbox\@tempboxa\par}%
   \else
      \hbox to\captionwidth{\hfil\box\@tempboxa\hfil}%
   \fi}%
\def\fignr{\small\sffamily\bfseries}%
\def\capt{\small\sffamily}%

\newdimen\@captionmargin\@captionmargin2cm\relax
\newdimen\captionwidth\captionwidth\hsize\relax

\let\reference\ref
\def\eqref#1{(\protect\ref{#1})}

\def\proof{\@ifnextchar[{\@proof}{\@proof[\unskip]}}
\def\@proof[#1]{\noindent{\bf Proof #1.}\enspace}

\def\myfrac#1#2{\raisebox{0.2em}{\small$#1$}\!/\!\raisebox{-0.2em}{\small$#2$}}
\def\abstractname{}

\@addtoreset {footnote}{page}

\renewcommand{\section}{%
   \@startsection
         {section}{1}{\z@}{-1.5ex \@plus -1ex \@minus -.2ex}%
               {1ex \@plus.2ex}{\large\bf}%
}
\renewcommand{\@seccntformat}[1]{\csname the#1\endcsname .
\quad}

{\let\@noitemerr\relax
\vskip-2.7em\kern0pt\begin{abstract}
\noindent{\bf Abstract.}\enspace
It is known that the minimal degree of the Jones polynomial 
of a positive knot is equal to its genus, 
and the minimal coefficient is $1$, with a similar relation for links.
We extend this result to almost positive links and partly identify the
3 following coefficients for special types of positive links.
We also give counterexamples to the
Jones polynomial-ribbon genus conjectures
for a quasipositive knot. Then we show that
the Alexander polynomial completely detects the minimal genus
and fiber property of canonical Seifert surfaces associated to
almost positive (and almost alternating) link diagrams.
\end{abstract}
}



\section{Introduction}

A link is called \em{quasipositive} if it is the closure of a braid
which is the product of conjugates of the Artin generators $\sg_i$
\cite{Rudolph2}. (We call such conjugates and their inverses positive
resp.\ negative \em{bands}.)
It is called \em{strongly quasipositive} if these conjugates are
positive \em{embedded bands} in the band representation of
\cite{Rudolph2}. It is called \em{positive} if it has a diagram with all
crossings positive (in the skein sense), and \em{braid positive} (or
a \em{positive braid link}) if it
has a braid representation which is a positive word in the
Artin generators. It is called \em{fibered}, if its complement in $S^3$
is a surface bundle over the circle.

Then
\begin{eqnarray}
\label{hie}\{\mbox{ quasipositive links }\}\,\supset\,
\{\mbox{ strongly quasipositive links }\}\,\supset\,
\{\mbox{ positive links }\}\,\supset\, & & \\[2mm]
\nonumber\rx{3cm}\,\supset\,\{\mbox{ fibered positive links }\}\,
\,\supset\,\{\mbox{ braid positive links }\}\,.
\end{eqnarray}
%
The only non-obvious inclusions are the second and fourth one. The
fourth inclusion is a well-known fact (it follows e.g. from
\cite{Gabai}), and the second inclusion
follows, as observed by Rudolph \cite{Rudolph} and Nakamura
\cite{Nakamura}, by applying the algorithm of Yamada
\cite{Yamada} or Vogel \cite{Vogel} to a positive diagram.
Links in some of the above classes have been studied, beside
by their intrinsic knot-theoretical interest, with
different motivations and in a variety of contexts,
including singularity theory \cite{ACampo,BoiWeb,Milnor},
algebraic curves \cite{Rudolph2,Rudolph3}, dynamical systems
\cite{BirWil} and (in some vague and yet-to-be understood way)
in 4-dimensional QFTs \cite{Kreimer}.

A different related class to positive links are the
almost positive links, those with almost positive
diagrams, which are, however, not positive. (A diagram
is almost positive if it has exactly one negative crossing.)

Let $g$ be the genus of a knot, $g_s$ the slice genus, and $g_r$
the ribbon genus. (Definitions are recalled below.)
For links similarly write $\chi$, $\chi_r$
and $\chi_s$ for the (Seifert), ribbon and slice Euler characteristic
resp. As any Seifert surface is a ribbon surface,
and any ribbon surface is (deformable into) a slice surface,
one has the inequalities $g\ge g_r\ge g_s$ and $\chi\le
\chi_r\le \chi_s$.

For knots we also have $u\ge g_r$, with $u$ being the \em{unknotting
number} \cite{Lickorish}. By the work of Kronheimer--Mrowka
\cite{KM,KM2} and Rudolph \cite{Rudolph2},
it is now known that the slice genus is estimated below by
the slice Bennequin inequality (a version of \cite[theorem 3]
{Bennequin} with $g$ replaced by $g_s$),
implying that for a strongly quasipositive knot $g=g_s$,
so that $u\ge g= g_r= g_s$. For positive braid knots $u\le g$
was known by \cite{BoiWeb}. Thus $u=g$ in this case.

Let $V$ be the Jones polynomial \cite{Jones}. Fiedler \cite{Fiedler}
proved that $\md V=g$ for a positive braid knot, and that
$\mcf V=1$. For positive braid links $L$ of $n=n(L)$ components,
$\mcf V=(-1)^{n-1}$ and $2\md V=1-\chi$. This follows
more generally for positive links $L$ by virtue of the fact that
positive diagrams are semiadequate (see \cite{LickThis}).
Fiedler further conjectured (his Conjecture 1)
that for arbitrary knots and links $L$, which
have a band representation on $s$ strands with $b$ bands,
\[
\md V\le \frac{b-s+1}{2}\,.
\]
He made a second conjecture (Conjecture 2), whose truth would imply
that equality in the above inequality is achieved only for
quasipositive links $L$.

In the paper of Kawamura \cite{Kawamura}, the theorems of
Fiedler and Kronheimer--Mrowka--Rudolph have been found to
imply that for a positive braid knot, $\md V=u$, with a
similar(ly obvious) relation for links. Then Kawamura quoted a special
case of Fiedler's first conjecture, asking whether it is true
(at least) for quasipositive links, and observing that the
slice Bennequin inequality would then imply the relation $\md V\le u$
for a quasipositive knot. (That $\md V=u$ does not extend to
quasipositive knots is easy to observe.)

In this paper we will give counterexamples to both Fiedler conjectures
of several special types, in particular the case of the first conjecture
addressed
by Kawamura. Then we will partly identify up to 3 of the coefficients
of the Jones polynomial of a positive link following the minimal one,
including a handy criterion to single out positive braid links, even
among fibered positive links. We will also extend Fiedler's result to
almost positive links. Some consequences are derived for the skein
polynomial \cite{HOMFLY}. Our results allow also to identify up to
two more coefficients of the skein polynomial of positive links.
For almost positive links, we obtain a proof of the inequality
conjectured by Morton \cite{Morton2} (for which in the general
case counterexamples are now known \cite{posex_bcr}).

We will use in some of the proofs the even valence graph version of the
Alexander polynomial studied in \cite{MS} with K.\ Murasugi. Applying
this method, we can also slightly improve and simplify the proof of
results of Hirasawa \cite{Hirasawa} and Goda--Hirasawa--Yamamoto
\cite{GHY}. These results are amplified by showing that the
Alexander polynomial completely determines the minimal genus
and fiber property of canonical Seifert surfaces associated
to almost positive (and almost alternating) link diagrams.
At the end we will give a few examples showing that many of
the possible extensions of these theorems are not true,
and mention some problems.

\section{Preliminaries}

\subsubsection*{Link polynomials}

The \em{skein polynomial} $P$ is a Laurent polynomial in two variables
$l$ and $m$ of oriented knots and links and can be defined
by being $1$ on the unknot and the (skein) relation
\begin{eqn}\label{1}
l^{-1}\,P\bigl(
\diag{5mm}{1}{1}{
\picmultivecline{-5 1 -1.0 0}{1 0}{0 1}
\picmultivecline{-5 1 -1.0 0}{0 0}{1 1}
}
\bigr)\,+\,
l \,P\bigl(
\diag{5mm}{1}{1}{
\picmultivecline{-5 1 -1 0}{0 0}{1 1}
\picmultivecline{-5 1 -1 0}{1 0}{0 1}
}
\bigr)\,=\,
-m\,P\bigl(
\diag{5mm}{1}{1}{
\piccirclevecarc{1.35 0.5}{0.7}{-230 -130}
\piccirclevecarc{-0.35 0.5}{0.7}{310 50}
}
\bigr)\,.
\end{eqn}
For a diagram $D$ of a link $L$, we will use all of the notations
$P(D)=P_D=P_D(l,m)=P(L)$ etc.\ for its skein polynomial, with
the self-suggestive meaning of indices and arguments.

The \em{Jones polynomial} $V$, and (one variable) \em{Alexander
polynomial} $\Dl$ are obtained from $P$ by the substitutions
\begin{eqnarray}
V(t) & = & P(-it,i(t^{-1/2}-t^{1/2}))\,, \label{PtoV}\\
\Delta(t) & \doteq & P(i,i(t^{1/2}-t^{-1/2})) \label{PtoDl}\,,
\end{eqnarray}
hence these polynomials also satisfy corresponding skein relations.
The sign `$\doteq$' means that the Alexander polynomial is
defined only up to units in $\bZ[t,t^{-1}]$; we will choose
the normalization depending on the context.

We denote in the following the coefficient of $t^m$ in $V(t)$ by $[V(t)]
_{m}$. In the case of a 2-variable polynomial, we index the bracket by
the whole monomial, and not just the power of the variables.
The \em{minimal} or \em{maximal degree} $\md V$ or $\Md V$ is the
minimal resp.\ maximal exponent of $t$ with non-zero coefficient in $V$.
An explicit (1-variable) polynomial may be denoted by the convention
of \cite{LickMil} by its coefficient list, when bracketing its absolute
term to indicate its minimal degree, e.g. $(-3\ [1]\ 2)\,=\,-3/t+1+2t$.
The \em{minimal} or \em{leading coefficient} $\mc V$ of $V$ is
$[V]_{\md V}$.

For an account on these link polynomials we refer to the
papers \cite{LickMil,Jones}. (Note: our convention here
for $P$ differs from \cite{LickMil} by interchange of
$l$ and $l^{-1}$, that is, our $P(l,m)$ is 
Lickorish and Millett's $P(l^{-1},m)$.)

\subsubsection*{Link diagrams}

A crossing $p$ in a link diagram $D$ is called \em{reducible}
(or \em{nugatory}) if $D$ can be represented in the form
\[
\diag{6mm}{4}{2}{
  \picrotate{-90}{\rbraid{-1 2}{1 1.4}}
  \picputtext[u]{2 0.7}{$p$}
  \picscale{1 1}{
    \picfilledcircle{0.7 1}{0.8}{$P$}
  }
  \picfilledcircle{3.3 1}{0.8}{$Q$}
}\es.
\]
$D$ is called reducible if it has a reducible crossing, else it is
called \em{reduced}.

A link diagram $D$ is \em{composite}, if there is a closed curve
$\gm$ intersecting (transversely) the curve of $D$ in two
points, such that both in- and exterior of $\gm$ contain crossings
of $D$, that is, $D$ has the form
\[
\diag{6mm}{4}{2}{
  \piccirclearc{2 0.5}{1.3}{45 135}
  \piccirclearc{2 1.5}{1.3}{-135 -45}
  \picfilledcircle{1 1}{0.8}{$A$}
  \picfilledcircle{3 1}{0.8}{$B$}
  {\piclinedash{0.1}{0.05}
   \piccircle{3.3 1}{1.3}{}
  }
}\quad. 
\]
Otherwise $D$ is \em{prime}. A link is prime if any in composite
diagram replacing one of $A$ and $B$ by a trivial (0-crossing)
arc gives an unknot diagram.

The diagram is \em{split}, if there is a closed curve not
intersecting it, but which contains parts of the diagram
in both its in- and exterior:
\[
\diag{6mm}{4}{2}{
  \picfilledcircle{1 1}{0.8}{$A$}
  \picfilledcircle{3 1}{0.8}{$B$}
  {\piclinedash{0.1}{0.05}
   \piccircle{3.3 1}{1.3}{}
  }
}\quad. 
\]
Otherwise $D$ is \em{connected} or \em{non-split}.
A link is split if it has a split diagram, and otherwise 
non-split.

We call a diagram $D$ \em{$k$-almost positive}, if
$D$ has exactly $k$ negative crossings. A link $L$ is
$k$-almost positive, if it has a $k$-almost positive diagram, but
no $l$-almost positive one for any $l<k$. We call a diagram or
link \em{positive}, if it is $0$-almost positive (see
\cite{Cromwell,Ozawa,Yokota,Zulli}), and \em{almost positive}
if it is $1$-almost positive \cite{apos}. Similarly one defines 
\em{$k$-almost negative}, and in particular \em{almost negative}
and \em{negative} links and diagrams to be the mirror images of
their $k$-almost positive (or almost positive or positive)
counterparts, and \em{($k$-)almost alternating} diagrams and links
\cite{Adams,Adams2}.

The \em{valency} of a Seifert circle $s$ is the number of
crossings attached to $s$. We call such crossings also
\em{adjacent} to $s$. 

\subsubsection*{Link surfaces}

A \em{Seifert} resp. \em{slice} surface of $L\subset S^3=
\partial B^4\subset B^4$
is a smoothly embedded compact orientable surface $S\subset S^3$
resp. $S\subset B^4$ with $\partial S=L$. A \em{ribbon surface}
is a smoothly immerged compact orientable surface $S\subset S^3$
with $\partial S=L$, embedded except at a finite number of double
transverse (ribbon) singularities, whose preimages are two arcs,
one lying entirely in the interior $\int S$ of $S$ and the other one,
too, except for its two endpoints, which lie on $\partial S$.
A \em{canonical (Seifert) surface} is a Seifert surface
obtained by Seifert's algorithm (see\cite{Rolfsen}).
We may allow (for links) all these surfaces to be disconnected, but 
they should have no closed ($\partial=\varnothing$) components.

The \em{(Seifert) genus} $g$, \em{slice genus} $g_s$, \em{canonical 
genus} $\tl g$ and \em{ribbon genus} $g_r$ are defined to be
the minimal genera of Seifert, slice, canonical resp.\ ribbon 
surfaces of $L$. Similarly one can define the (Seifert), slice,
canonical resp.\ ribbon \em{Euler characteristic} $\chi$, $\chi_s$,
$\tl\chi$, and $\chi_r$ to be the maximal Euler characteristic
of such surfaces of $L$.

In \cite[theorem 3]{Bennequin}, Bennequin shows that
for a braid $\bt$ on $s(\bt)$ strands, with writhe (exponent sum)
$w(\bt)$ and with closure $\hat\bt=K$, we have an estimate for the 
Euler characteristic $\chi(K)$ of $K$:
\[
1-\,\chi(K)\,\ge\,w(\bt)\,-\,s(\bt)\,+1\,.
\]
This is easily observed to extend by means of the algorithm of Yamada
\cite{Yamada} or Vogel \cite{Vogel} to an inequality for arbitrary
link diagrams $D$ of $K$:
\begin{eqn}\label{b}
1-\,\chi(K)\,\ge\,w(D)\,-\,s(D)\,+1\,=:\,b(D)\,,
\end{eqn}
with $w(D)$ being the writhe of $D$, and $s(D)$ the number of its
Seifert circles. We call the r.h.s. of \eqref{b} the \em{Bennequin
number} of $D$. (It clearly depends a lot on the diagram for
a given link.)

Rudolph \cite{Rudolph3} later improved this inequality, by replacing
$\chi(K)$ by $\chi_s(K)$: 
\begin{eqn}\label{b2}
1-\,\chi_s(K)\,\ge\,b(D)\,.
\end{eqn}

Recently, he obtained a further improvement, this time by raising the
r.h.s. \cite{Rudolph2}:
\begin{eqn}\label{rb}
1-\,\chi_s(K)\,\ge\,w(D)\,-\,s(D)\,+1\,+2s_-(D)=:\,rb(D)\,,
\end{eqn}
with $s_-(D)$ being the number of ($\ge 2$-valent) Seifert circles of
$D$, to which only negative non-nugatory crossings are adjacent. The
new quantity on the right we call
\em{Rudolph-Bennequin number} of $D$. Again $rb(D)$
heavily depends on the diagram, even more than $b(D)$.
(For example, unlike $b(D)$, $rb(D)$ is no longer invariant under
flypes and mutations.) Thus again one is interested in choosing
for a given link $K$ the diagram $D$ so that $rb(D)$ is as
large as possible.

\section{Counterexamples to the Jones polynomial-ribbon genus
conjectures}

\subsection{Preparations}

While the improvement \eqref{rb}, as compared to \eqref{b2}, may
not seem significant at first sight, it has the advantage of
eliminating the minimal $l$-degree in the skein polynomial
$\md_lP$ as an obstruction to increasing the estimate by proper
choice of the diagram $D$, since by \cite{Morton} we always have
$b(D)\le \md_lP(K)$. 

A practical example where this turned out helpful was given in
\cite{4gen}, and is recalled below, as it will be used.
(The notation for knots we apply is the one of Rolfsen's tables
\cite[appendix]{Rolfsen} for $\le 10$ crossings, and of
the knot table program KnotScape \cite{HT} for $11$ to $16$ 
crossings. By $!K$ we will denote the obverse, or mirror image, of $K$.)

\begin{exam}
The knot $13_{6374}$ has $\md_lP=0$ and Alexander polynomial 
$\Dl=1$. It has many
diagrams $D$ with $b(D)=0$, but it cannot have any such diagram with
$b(D)>0$, because of Morton's inequality. However, it does have
diagrams $D$ with $rb(D)>0$, thus showing it not to be slice.
\end{exam}

\begin{figure}[htb]
\[
\begin{array}{c} \epsfsv{4cm}{t1-13_6374} \\ \ry{6mm}13_{6374}
\end{array}
\]
\caption{A non-slice knot with unit Alexander polynomial.\label{ppf}}
\end{figure}

In order to construct our counterexamples, we need a few more
simple lemmas.

\begin{lemma}
\label{lemm1}
If $K$ is quasipositive, then $\md_lP(K)\ge 1-\chi_s(K)$.
\end{lemma}

\proof If $D$ is a diagram of a quasipositive braid representation
of $K$, then $1-\chi_s(K)=b(D)$, and $b(D)\le \md_lP(K)$ by Morton's
inequality. \qed

In the following $K_1\# K_2$ denotes the connected sum of
$K_1$ and $K_2$, and $\#^nK$ denotes the connected sum of $n$
copies of $K$.

\begin{lemma}
\label{lemm2}
If $K_{1,2}$ have diagrams $D_{1,2}$ which are not negative, then
$K_1\# K_2$ has a diagram $D$ with $rb(D)=rb(D_1)+rb(D_2)$.
\end{lemma}

\proof We must apply the connected sum of $D_{1,2}$ so that the
Seifert circles of $D_{1,2}$ affected by the operation have at
least one positive crossing adjacent to them. \qed

\begin{lemma}
\label{lemm3}
If $K$ is strongly quasipositive, then $\chi(K)=\chi_s(K)$.
\end{lemma}

\proof For the Seifert surface $S$ associated to
a strongly quasipositive braid representation diagram $D$ of $K$,
we have
\[
1-\chi(K)\le 1-\chi(S)\,=\,b(D)\,\le\,rb(D)\,\le\,1-\chi_s(K)\,\le\,
1-\chi(K)\,,
\]
implying equality everywhere. \qed

\subsection{Degree inequality conjecture}

Fiedler's first conjecture was whether
\[
\md V_L\,\le\, \frac{b-s+1}{2}\,,
\]
if $L$ has a $b$-band representation on $s$ strands, and
Kawamura's (weaker) question was whether it is true at least
if this band representation is positive.

\begin{figure}[htb]
\[
\begin{array}{c}
\epsfsv{4.6cm}{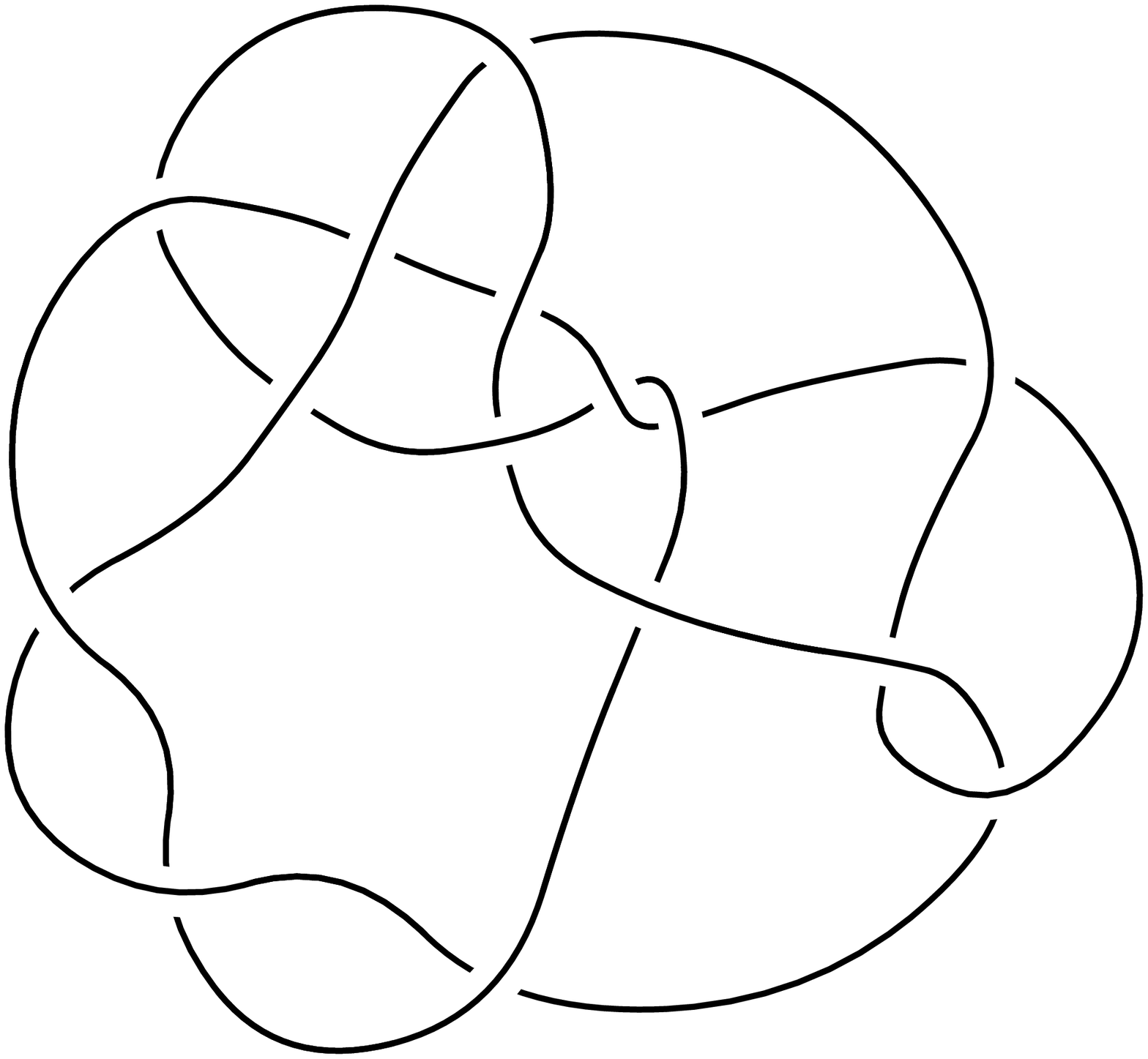} \\[2.2cm]
!15_{162508}
\end{array}
\]
\caption{\label{fig1}}
\end{figure}

\begin{exam}
Consider the knot $!15_{162508}$ (see figure \reference{fig1}).
Using the method described in \cite[appendix]{deg}, it was
found that it is ribbon (and hence slice), and one calculates $\md V=1$.
It turns out to have the quasipositive 5-braid representation
\[
(\sg_1^{-1}\sg_2^{-1}\sg_3\sg_4\sg_3^{-1}\sg_2\sg_1)\,
(\sg_2^{-1}\sg_1\sg_2)\,(\sg_2\sg_3^{-1}\sg_4\sg_3\sg_2^{-1})\,\sg_3\,.
\]
(The knot can be identified from this representation
by the tool {\tt knotfind} included in \cite{HT}. Note that
this representation also directly shows sliceness.)
Thus it is a slice example answering negatively Kawamura's
question, and hence also a counterexample Fiedler's first conjecture.
\end{exam}

Another special type of example is 

\begin{figure}[htb]
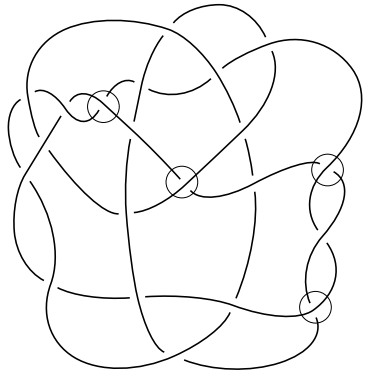

\[
\diag{1cm}{5.0}{5.0}{
  \picputtext{2.5 2.5}{\epsfs{5cm}{t1-morton1}}
  \piclinewidth{100}
  \piccircle{2.32 2.65}{0.21}{}
  \piccircle{1.28 3.65}{0.21}{}
  \piccircle{4.25 2.81}{0.21}{}
  \piccircle{4.09 0.99}{0.21}{}
}
\]
\caption{The knot $K$ in example \reference{ex2}.\label{fig2}}
\end{figure}

\begin{exam}\label{ex2}
Consider the knot $K$ in figure \reference{fig2},
which is the closure of the $4$-braid
\begin{eqn}\label{brep}
\sg_1^{2}(\sg_1\sg_2\sg_1^{-1})\sg_2\sg_1\sg_3
(\sg_1\sg_2\sg_1^{-1})\sg_2(\sg_2\sg_3\sg_2^{-1})
(\sg_1\sg_2\sg_1^{-1})(\sg_2\sg_3\sg_2^{-1})\,.
\end{eqn}
This braid is quasipositive, in fact, strongly quasipositive.
The diagram of $K$ in figure \reference{fig2} was obtained
from that representation. 

One easily sees that $g=g_s=4$. But $\md V=5$. Thus $g_s<\md V$.
In fact, this knot has unknotting number $4$. (Switch the
encircled crossings in the diagram of figure \reference{fig2}.)
Thus even the weaker inequality, in which Kawamura was
interested, $\md V\le u$ is not always true.
\end{exam}

\begin{rem}
$!15_{162508}$ is surely not strongly quasipositive, as $g>0=g_s$.
Thus the above example $K$ is the most special in the hierarchy
\eqref{hie}.
\end{rem}

The only case of some interest, remaining not covered by the
above examples, is that of a slice knot with $u<\md V$. Very
likely such examples exist, too, although I didn't find any.

\begin{rem}
If one is interested in a \em{general} knot $K$ with $\md V>u$, then
there is a much simpler and well-known example, $!10_{132}$. It has
$u=1$, but $\md V=2$. However, $!10_{132}$ is not quasipositive.
As it is not ribbon, or slice (its determinant 5 is not a
square), it has 4-genus $1$, and a quasipositive representation
of $n$ strands would have $n+1$ bands. Then the untwisted
2-cable link $(!10_{132})_2$ would have a representation
on $2n$ strands of writhe $2n+2$. Thus by \cite{Morton},
$\md_lP((!10_{132})_2)\ge 3$, but from the calculation of
\cite{MortonShort} we know $\md_lP((!10_{132})_2)=1$.
\end{rem}

\begin{rem}
T.\ Tanaka, in a preprint \cite{Tanaka} sent to T.\ Fiedler, has
claimed independently counterexamples to Fiedler's first conjecture.
On the opposite end, M.\ Ishikawa \cite{Ishikawa} proved Fiedler's
inequality for some links obtained by A'Campo's method \cite{ACampo}.
\end{rem}

\subsection{Extremal property conjecture}

Fiedler also conjectured (his Conjecture 2) that if
a link $L$ has a $b$-band $s$-strand band representation
with
\begin{eqn}\label{c1}
\md V_L\,=\,\frac{b-s+1}{2}\,,
\end{eqn}
then it is quasi-positive. (Fiedler's formulation is slightly
different, but easily implies the one given here.)

We will now construct a counterexample also to this conjecture,
albeit some more effort is necessary, and we must use the
example found previously in a related context in \cite{4gen}.
Our counterexample has likely crossing number 58.

\begin{prop}\label{ppp1}
The knot $K'=13_{6374}\,\#\,\#^3(!15_{162508})$,
is not quasipositive, yet it has a band representation with
equality in \eqref{c1}.
\end{prop}

\proof We first discuss the prime factors separately.

\def\labelenumi{\arabic{enumi}.}
\def\theenumi{\arabic{enumi}.}
\begin{enumerate}

\item Consider $13_{6374}$. By switching one of the crossings
in the clasp in the lower right part of the diagram in figure
\reference{ppf}, one obtains
$4_1$. Thus $\#^213_{6374}$ turns by 2 crossing changes into
the slice knot $4_1\#4_1$. Hence $1-\chi_s(\#^2 13_{6374})\le 4$.
On the other hand, as $13_{6374}$ has a diagram $D$ with $rb(D)=2$,
we have by lemma \reference{lemm2} that $1-\chi_s(\#^2 13_{6374})=4$.

As said, 
\[
\md_lP(\#^2 13_{6374})\,=\,0\,=\frac{1-\chi_s(\#^2 13_{6374})}{2}-2\,.
\]
Also $\md V=-1$ by calculation. As $\Md_mP(13_{6374})=4$ and it has
crossing number $<15$, by \cite{gen2} $13_{6374}$ has a diagram of
canonical genus $2$, and thus by applying Yamada's algorithm
\cite{Yamada} on it, we obtain a(n embedded) band representation with
\[
\frac{b-s+1}{2}=2\,=\md V+3\,.
\]

\item For $!15_{162508}$, we have a quasipositive band representation
as $5$-braid with $4$ bands, and it is slice. Thus
\[
1-\chi_s\,=\,0\,=\,\frac{\md_lP}{2}\quad\mbox{and}
\quad \frac{b-s+1}{2}\,=\,0\,=\,\md V-1\,.
\]
\end{enumerate}

In summary we have the following situation for proper diagrams $D$
and $b$-band $s$-braid representations:
\[
\begin{mytab}{c||c|c}{&&}
\ry{1.5em} & $13_{6374}$ & $!15_{162508}$ \\[2mm]
\hline
\hline
\ry{1.5em} $\Myfrac{\md_lP}{2}-\Myfrac{rb(D)}{2}$ & $-1$ & $0$ \\[2mm]
\hline
\ry{1.5em} $\md V\,-\,\Myfrac{b-s+1}{2}$ & $-3$ & $1$ \\[2mm]
\end{mytab}
\]

Since both quantities are additive under connected sum for
proper diagrams and band representations
(by Lemma \reference{lemm1} resp. in the obvious way),
we obtain for $K'$ a band representation
with 
\[
\frac{b-s+1}{2}\,=\,\md V(K')\,,
\]
but also a diagram $D$ of $K'\# K'$ with
\[
\md_l P(K'\# K')\,<\,rb(D)\,\le\,1-\chi_s(K'\# K')\,,
\]
(in fact $rb(D)=1-\chi_s(K'\# K')$), so that $K'\# K'$ is not
quasipositive by lemma \reference{lemm1}. Then $K'$ cannot
be quasipositive either. \qed

%
%

One can also obtain a counterexample to an ``embedded band''
version of Fiedler's conjecture, namely whether a knot $K$
with an \em{embedded} band representation achieving equality in
\eqref{c1} is \em{strongly} quasipositive.

\begin{prop}
The knot $K'=\#^3K\,\#\,6_1$,
with $K$ being the knot in example \reference{ex2},
is not strongly quasipositive, yet it has an embedded
band representation satisfying equality in \eqref{c1}.
\end{prop}

\proof As $6_1$ has canonical genus $1$, it has an embedded
band representation with
\[
\frac{b-s+1}{2}\,=\,1\,=\,\md V\,+3\,.
\]
Now consider $K$. It has a strongly quasipositive 
band representation with $b=11$ bands on $s=4$ strands, so
that
\[
1-\chi_s\,=\,\frac{b-s+1}{2}\,=\,4\,.
\]
However, 
\[
\md V\,=\,5\,.
\]
Thus $K'$ has an embedded 
band representation satisfying equality in \eqref{c1}.

As genus is additive under connected sum, we have $g(K')=13$.
However, as $g_s$ is subadditive under connected sum, and $6_1$
is slice, we have $g_s(K')\le 12$, so that $g>g_s$, and so $K'$
is not strongly quasipositive by lemma \reference{lemm3}. \qed


%
%

There is an exponentiated version of Fiedler's conjecture,
namely to ask about (non-strong) quasipositivity assuming
equality in \eqref{c1} for an embedded band representation.
We conclude this section by showing how to construct
counterexamples also for this most sharp case.

The problem to give such counterexamples reduces in one possible way
to replacing $!15_{162508}$ by a strongly quasipositive knot with
$\md V>g$. Then the same argument as in the proof of proposition
\reference{ppp1} goes through with embedded band representations.

\begin{figure}[htb]
\[
\begin{array}{c}
\epsfsv{4.6cm}{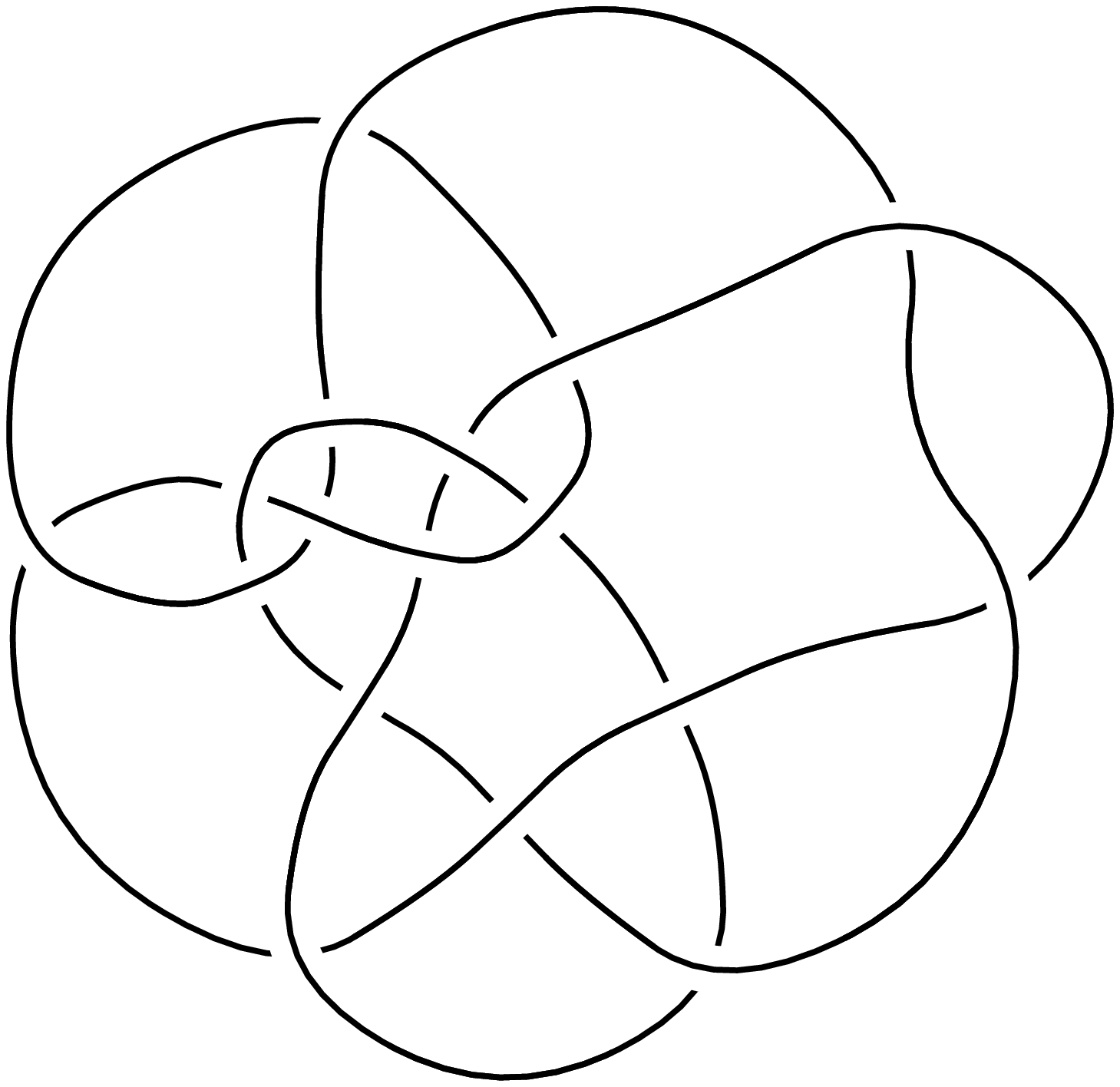} \\[2.0cm]
\end{array}
\]
\caption{\label{fig_3}}
\end{figure}

\begin{exam}\label{Cx}
Consider the (apparently) 17 crossing knot on figure \reference{fig_3}.
It has a band representation with 7 bands on 4 strands,
\[
\bigl(\,(\sg_2\sg_3\sg_2^{-1})(\sg_1\sg_2\sg_1^{-1})\,\bigr)^3\sg_1.
\]
(The diagram in figure \reference{fig_3} was obtained again using
KnotScape.) Thus its genus is $g=2$. Also $\md_lP=4$, but $\md V=3$.
\end{exam}

\begin{rem}
It is clear from example \reference{Cx} that in fact we could have
used it also as counterexample to Fiedler's first conjecture.
However, unlike for the knot $K$ in example \reference{ex2},
I cannot show $u=g(=2)$ here. On the other hand, $K$ cannot
be used in example \reference{Cx}, because it has $\md_lP=10$.
($K$ was found as a counterexample to Morton's conjecture,
as reported in \cite{posex_bcr}.) This way, any of the previous knots
comes to its own right.
\end{rem}

\section{The coefficients of the Jones polynomial}

\noindent{\bf Convention.}
It is convenient to assume in the sequel that all diagrams
we consider are non-split. In particular, since
non-split positive diagrams represent non-split links, we assume
all positive links to be so.

\begin{defi} 
A \em{separating} Seifert circle is a Seifert circle 
with non-empty interior and exterior. (That is, both
interior and exterior contain crossings, or equivalently,
other Seifert circles.) A diagram with no
separating Seifert circles is called \em{special}. 
\end{defi}

Any diagram decomposes as \em{Murasugi sum} along its separating
Seifert circles into special diagrams (see \cite[\S 1]{Cromwell}).
{}For any diagram, any two of the properties positive, alternating and
special imply the third. We call these diagrams \em{special 
alternating} to conform to the classical terminology of Murasugi
\cite{Murasugi3}.

\subsection{Positive braids}

It was known (e.g. from \cite{Fiedler}) that the minimal
coefficient of the Jones polynomial of a positive braid link
is $\pm 1$. We will show here a statement on the next 3 coefficients.

\begin{theo}\label{th1}
Let $L$ be a non-split braid positive link of $n(L)$ components.
Then
\[
(-1)^{n(L)-1}t^{(\chi(L)-1)/2}V_L(t)\,=\,1+pt^2+kt^3+(\mbox{ higher
order terms })\,,
\]
with $p=p(L)$ being the number of prime factors of $L$ and 
\begin{eqn}\label{k}
-p\es\le\es k\es\le\es \frac 32\,\bigl( 1-\chi(L)-p \bigr)\,,
\end{eqn}
where $\chi(L)$ is the Euler characteristic of $L$.
\end{theo}

\begin{figure}[htb]
\[
\begin{array}{*4c}
\epsfsv{3.8cm}{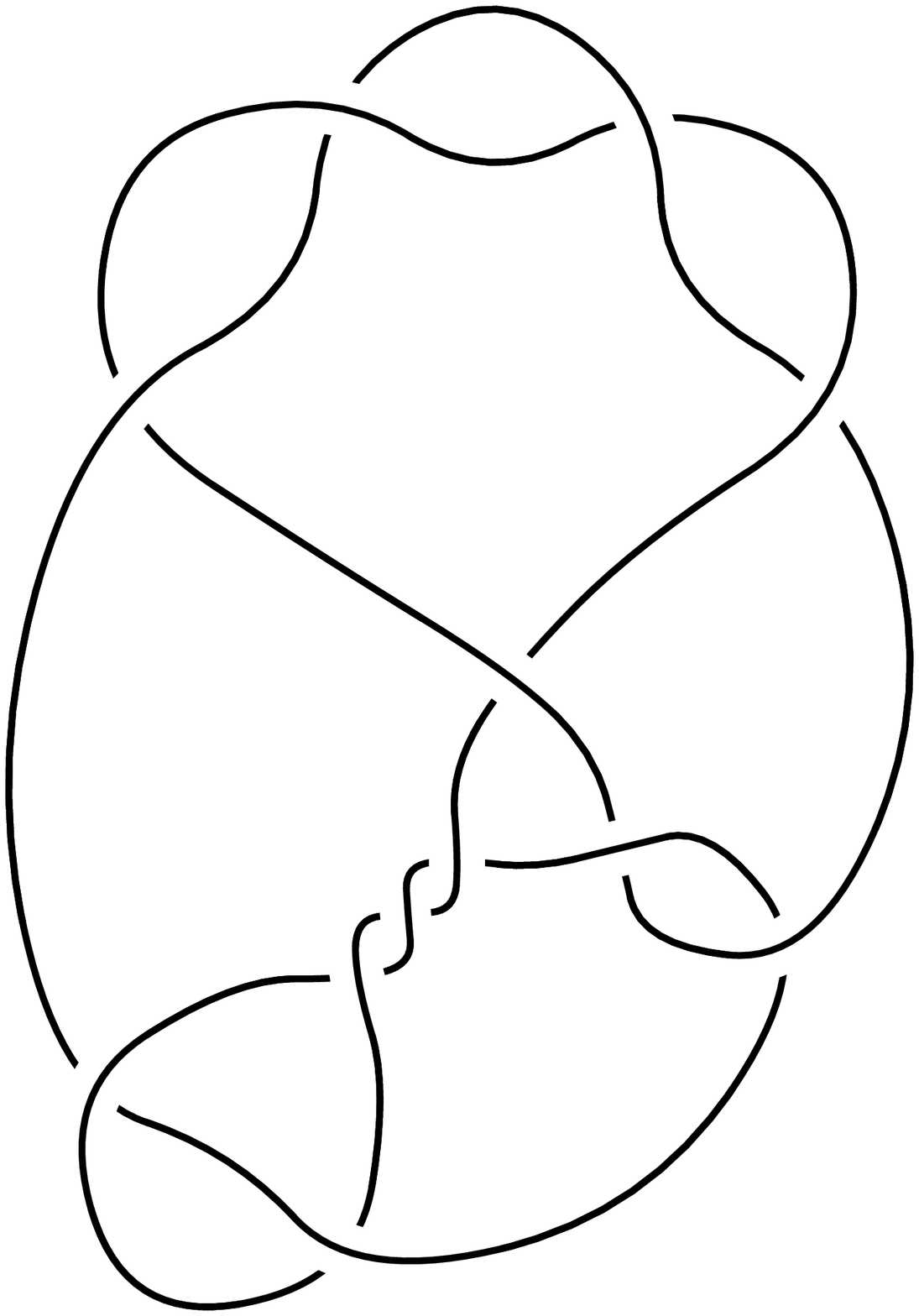} & \epsfsv{3.8cm}{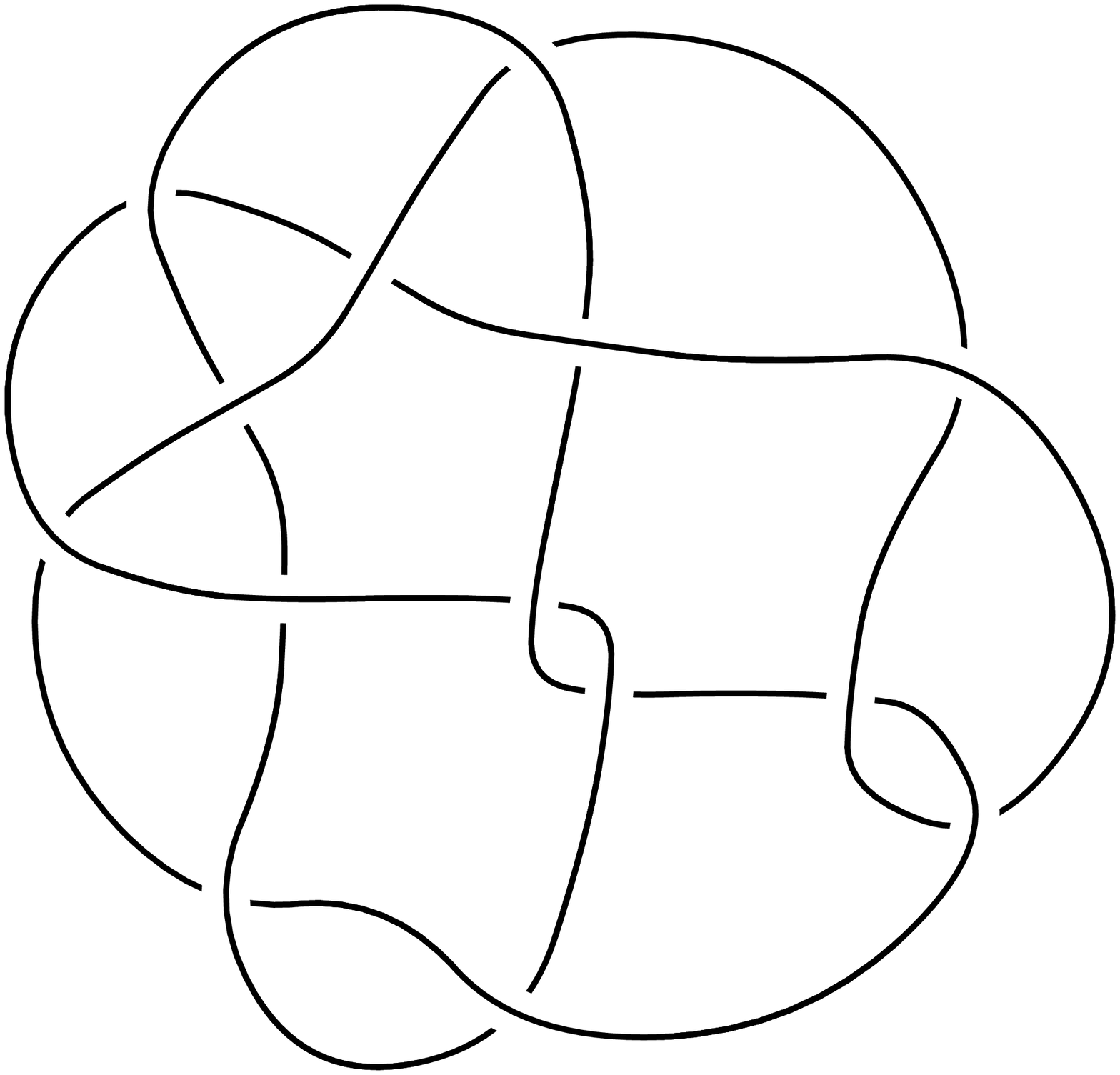} &
\epsfsv{3.8cm}{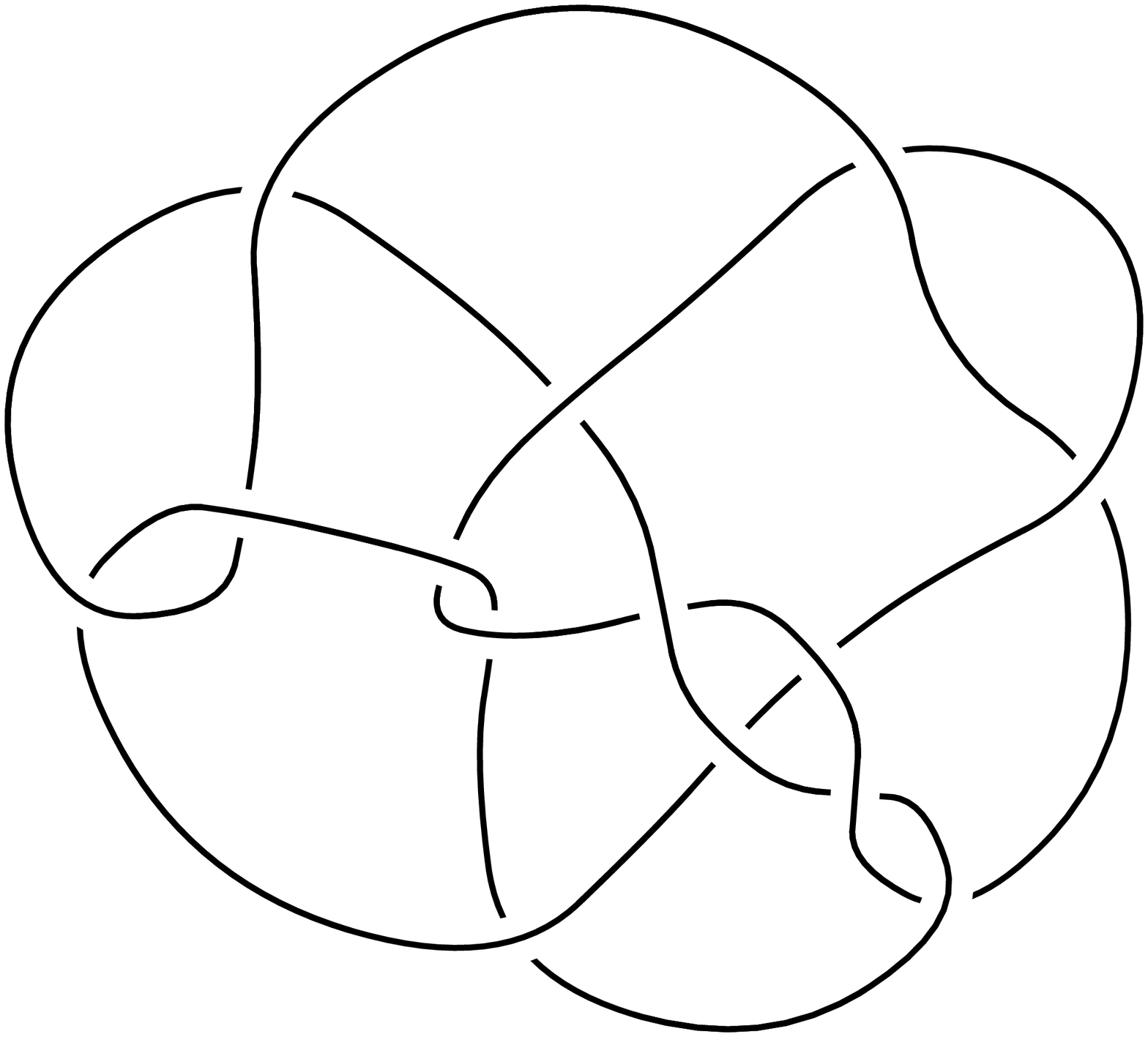} & \epsfsv{3.8cm}{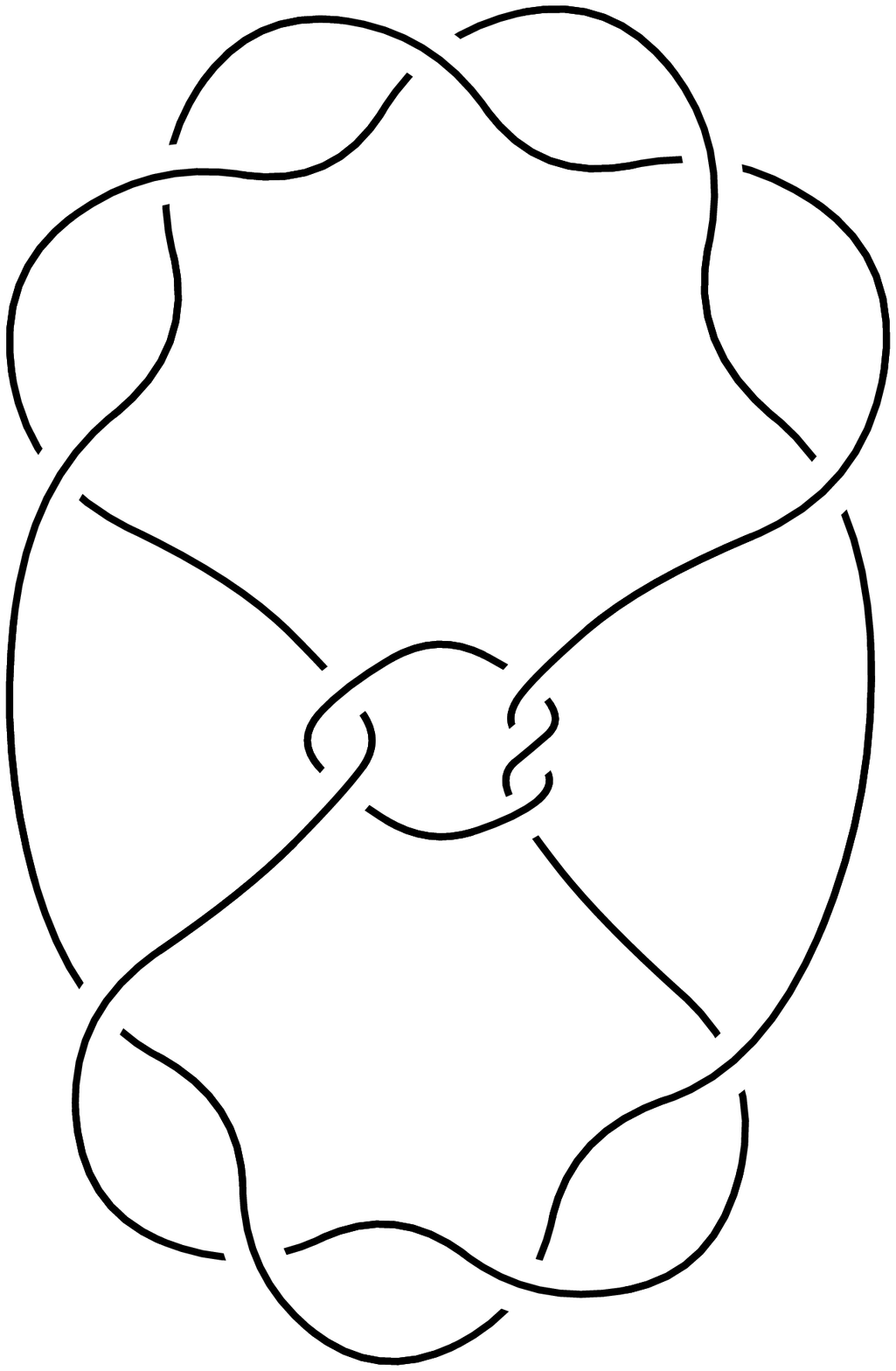} \\[2.2cm]
12_{420} & 14_{33805} & 14_{37899} & 14_{39977}
\end{array}
\]
\caption{\label{fig3}}
\end{figure}

Note that it is a rather unusual situation to be able to
read the number of prime factors off the polynomial.
This is, for example, not possible for alternating links
as shows the well-known pair $8_9$ and $4_1\#4_1$~--
the one knot is prime and the other one composite, yet they have
the same Jones polynomial. Two more interesting examples of
this type are as follows:

\begin{exam}
With some effort one also finds such pairs of positive (or special)
alternating knots: $12_{420}$ (figure \reference{fig3}) and
$!3_1\#9_{13}$ or $14_{4132}$ and $!5_2\#!9_9$.
\end{exam}

\begin{exam}
Even more complicated, but still existent, are such examples
of fibered positive knots. The simplest group I found
is a triple consisting of $14_{39977}$, $!3_1\#14_{33805}$ and
$!3_1\#14_{37899}$ (see figure \reference{fig3}).
\end{exam}

\proof[of theorem \reference{th1}]
If $\bt$ is a positive braid diagram of $L$, then by the result of
\cite{Cromwell2} the number of prime factors of the link $\hat\bt$
is equal to the number of prime factors of the diagram $\hat\bt$.
By \cite{BoiWeb} one can always choose $\bt$ so that it contains
a $\sg_i^2$. 
Apply the skein relation at one of the crossings. Then
\[
V_+\,=\,t^2V_-\,+\,(t^{3/2}-t^{1/2})V_0\,,
\]
with $L_-$ and $L_0$ both braid positive. Let $p_*=p(L_*)$.
By induction on the crossing number of the braid we have
\begin{eqn}\label{sum}
\begin{array}{*8c}
t^2 V_- & = & (-1)^{n(L)-1}t^{(1-\chi(L))/2}\,\cdot &
(\es[0]\ & 1\ & 0 \ & p_- \ & \dots\es) \\
(t^{3/2}-t^{1/2})V_0
& = & (-1)^{n(L)-1}t^{(1-\chi(L))/2}\,\cdot 
& (\es[1]\ & -1\ & p_0\ & k_0-p_0\ & \dots\es) \\
& & & \multicolumn{5}{c}{\rule[0.04cm]{5.8cm}{0.02cm}} \\
V_+ & = & (-1)^{n(L)-1}t^{(1-\chi(L))/2}\,\cdot &
(\es[1]\ & 0\ & p_0\ & k_0-p_0+p_-\ & \dots\es)
\end{array}
\end{eqn}
As $p_+=p_0$ and $0\le p_--p_+=p_--p_0\le 2$, the claim
follows by induction, once it is checked directly for
connected sums of trefoils and Hopf links, except for the
right inequality in  \eqref{k}, which follows only with
the constant $\myfrac{3}{2}$ replaced by $2$. (Note, that
$k$ and $p$ are both additive under connected sum.)
To prove $k\le \myfrac{3}{2}(1-\chi-p)$, we need to show
that after a smoothing with $p_-=p_++2$ we can choose
another one with $p_-\le p_++1$.

Write 
\[
\bt=\ds\prod_{k=1}^l\sg_i^{m_k}w_k
\]
with all $w_k$ containing
no $\sg_i$ but some of $\sg_{i\pm 1}$. Then one of the $k_i$, say
$k_1$, is equal to $2$, $k_2\ge 2$ and $l=2$. Then after
smoothing out one of the crossings in the clasp, we have $k_1=1$,
and then applying the skein relation at the other clasp, we
have $p_--p_0\le 1$, as desired. \qed

{}From the proof it is clear that the second inequality in \eqref{k}
is not sharp, and with some work it may be improvable. Candidates
for the highest ratio $k/(1-\chi-p)$ are braids of the form
$(\sg_1^2\sg_2^2\dots\sg_l^2)^2$, for which with $l\to\infty$
this ratio converges upward to $1$.

Contrarily, the first inequality is clearly sharp, namely for
connected sums of $(2,.)$-torus links.

\begin{question}\label{qu1}
Are the only links with $p+k=0$ connected sums of $(2,.)$-torus links?
\end{question}

\subsection{Fibered positive links}

We shall now prove a result on almost positive diagrams, which
shows a weaker version of theorem \reference{th1} for fibered positive
links. We need one definition.

\begin{defi}
The \em{Seifert graph} $\bar S_D$ of a diagram $D$ is a graph
obtained by putting a vertex for each Seifert circle of $D$
and connecting by an edge two vertices if a crossing is
joining the two corresponding Seifert circles. (If two
Seifert circles are connected by several crossings, $\bar S_D$
has multiple edges.) The \em{reduced Seifert graph} $S_D$
of $D$ is obtained by removing edges of $\bar S_D$ such that
(a) $S_D$ has no multiple edge and (b) two vertices are
connected by an edge in $\bar S_D$ iff they are so in $S_D$.
\end{defi}

\begin{defi}
For a link diagram $D$, let
$\chi(D)=s(D)-c(D)$, where $s(D)$ is the number of
\em{Seifert circles}, and $c(D)$ the number of \em{crossings} of $D$.
$\chi(D)$ is the \em{canonical Euler characteristic} of $D$.
\end{defi}

\begin{theorem}\label{tht}
Let $D$ be an almost positive diagram of a link $L$ with $n(L)$
components, with negative crossing $p$. If there is another
crossing in $D$ joining the same two Seifert circles as $p$,
then $\md V_L\ge (1-\chi(D))/2$. Otherwise, $\md V_L=(1-\chi(D))/2-1$
and $\mcf V_L=(-1)^{n(L)-1}$.
\end{theorem}

\begin{figure}[htb]
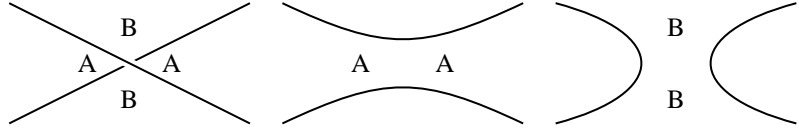

\[
\diag{8mm}{4}{2}{
   \picline{0 0}{4 2}
   \picmultiline{-5.0 1 -1.0 0}{0 2}{4 0}
   \picputtext{2.7 1}{A}
   \picputtext{1.3 1}{A}
   \picputtext{2 1.6}{B}
   \picputtext{2 0.4}{B}
} \quad
\diag{8mm}{4}{2}{
   \pictranslate{2 1}{
       \picmultigraphics[S]{2}{1 -1}{
           \piccurve{-2 1}{-0.3 0.2}{0.3 0.2}{2 1}
       }
   }
   \picputtext{2.7 1}{A}
   \picputtext{1.3 1}{A}
} \quad
\diag{8mm}{4}{2}{
   \pictranslate{2 1}{
       \picmultigraphics[S]{2}{-1 1}{
           \piccurve{2 -1}{0.1 -0.5}{0.1 0.5}{2 1}
       }
   }
   \picputtext{2 1.6}{B}
   \picputtext{2 0.4}{B}
}
\]
\caption{\label{figsplit}The A- and B-corners of a
crossing, and its both splittings. The corner A (resp. B)
is the one passed by the overcrossing strand when rotated 
counterclockwise (resp. clockwise) towards the undercrossing 
strand. A type A (resp.\ B) splitting is obtained by connecting 
the A (resp.\ B) corners of the crossing.}
\end{figure}

Recall, that the Kauffman bracket $[D]$ \cite{Kauffman} of a
link diagram $D$ is a Laurent polynomial in a variable
$A$, obtained by summing over all states $S$ the terms
\begin{eqn}\label{eq_12}
A^{\#A(S)-\#B(S)}\,\left(-A^2-A^{-2}\right)^{|S|-1}\,,
\end{eqn}
where a state is a choice of splittings of type $A$ or 
$B$ for any single crossing (see figure \ref{figsplit}), 
$\#A(S)$ and $\#B(S)$ denote the number of
type A (resp. type B) splittings and $|S|$ the
number of (disjoint) circles obtained after all
splittings in a state.

The Jones polynomial of a link $L$ is related to the
Kauffman bracket of some diagram of it $D$ by
\begin{eqn}\label{conv}
V_L(t)\,=\,\left(-t^{-3/4}\right)^{-w(D)}\,[D]
\raisebox{-0.6em}{$\Big |_{A=t^{-1/4}}$}\,,
\end{eqn}
$w(D)$ being the writhe of $D$.

\proof The maximal possible degree of $A$ in 
\begin{eqn}\label{S}
[D]\,=\,\sum_{S\scbox{ state}}
A^{\#A(S)-\#B(S)}\,\left(-A^2-A^{-2}\right)^{|S|-1}\,,
\end{eqn}
is that of the $A$-state (the state with all crossings
$A$-splitted), as under any splitting switch $A\to B$,
the power of $A$ in first factor in \eqref{eq_12} goes down by $2$,
and the maximal power of $A$ in second factor in \eqref{eq_12}
increases at most by $2$. If $D$ is almost positive with
negative crossing $p$, then the maximal possible power of $A$
in \eqref{S} is $A^{c(D)+2(s(D)-2)}$, as the $A$-state $S_A$ has
$s(D)-1$ loops.  They are the Seifert circles not adjacent to
$p$, and a loop consisting of the 2 Seifert circles, call
them $a$ and $b$, adjacent to $p$.

Now we must consider what states contribute terms of
$A^{c(D)+2(s(D)-2)}$ in \eqref{S}. These are exactly the states,
for which, when obtained from the $A$-state by successively switching
$A\to B$ splittings, $|\,.\,|$ increases under any such switch.

Let $\br{S\,:\,k}\in\{A,B\}$ be the split of $k$ in $S$, and
let $s_k(S)$ be the state obtained by switching splitting $A\to B$
at crossing $k$ in $S$, assuming $\br{S\,:\,k}=A$. Then if
$|s_k(S_A)|<|S_A|$, any state $S$ with $\br{S\,:\,k}=B$ is
irrelevant for the highest term in \eqref{S}. Clearly, this
happens whenever $k$ is a crossing connecting one or two
Seifert circles not adjacent to $p$. 

Thus the only terms contributing to $A^{c(D)+2(s(D)-2)}$ in \eqref{S}
are those for which $\br{S\,:\,k}=B$ implies that $k$
has the same two adjacent Seifert circles $a$ and $b$ as $p$.

Let $p_1,\dots,p_k=p$ be these crossings. Since any splitting
switch $A\to B$ in $s_p(S_A)$ reduces $|\,.\,|$, the only
state $S$ with $\br{S\,:\,p}=B$ relevant for the highest term in
\eqref{S} is $s_p(S_A)$, whose contribution to the coefficient
of this highest term is $(-1)^{|s_p(S_A)|-1}=(-1)^{s(D)-1}$.

It is also easy to see that if $\br{S\,:\,p}=A$, any of the
$2^{k-1}$ remaining states $S$ to consider contribute to
$A^{c(D)+2(s(D)-2)}$, the coefficient being $(-1)^{s(D)+\#B(S)}$,
as $|S|=s(D)-1+\#B(S)$.
The sum over all such $S$ of these coefficients is $(-1)^{s(D)}$
times the alternating sum of binomial coefficients. Thus this sum
vanishes for $k-1>0$, and cancels for $k-1=0$ the coefficient
$(-1)^{s(D)-1}$ of $s_p(S_A)$.

The rest follows from \eqref{conv} with $w(D)=c(D)-2$, and
the remark that $1-\chi(D)$ and $n(L)-1$ have the same parity. \qed

\begin{corr}\label{corr1}
Let $L$ be a fibered positive link of $n(L)$ components.
Then $[V_L(t)]_{(3-\chi(L))/2}=0$, that is,
\[
(-1)^{n(L)-1}t^{(\chi(L)-1)/2}V_L(t)\,=\,1+kt^2+(\mbox{ higher
order terms })\,,
\]
with $k$ being some integer.
\end{corr}

\proof This is proved as theorem \reference{th1}
by induction on the crossing number of a positive diagram $D$.
Apply the skein relation at any (non-nugatory) crossing $p$ of $D$.
Since the reduced Seifert graph of $D$ is a tree, there is
another crossing between the same two Seifert circles. Let $D_0$ be
$D$ with $p$ smoothed out, and $L_0$ be the link $D_0$ represents.
$L_0$ is still fibered, because $D_0$ is positive and connected,
and its reduced Seifert graph is still a tree. Similarly let $D_-$ be
$D$ with $p$ switched and $D_-$ representing a link $L_-$.
Then to $L_-$ we can apply the above theorem. So
$\md V_-=(1-\chi)/2-1$, and the coefficients of
$t^{(1-\chi)/2+1}$ in $t^2V_-$ and $(t^{3/2}-t^{1/2})V_0$
cancel as in \eqref{sum}. \qed

\subsection{Positive and almost positive links}

Corollary \reference{corr1} is a special
case of the following result, describing the second coefficient
of the Jones polynomial for an arbitrary positive link.

\begin{theorem}\label{cc}
Let $L$ be a positive link with positive diagram $D$. Then 
\[
(-1)^{n(L)-1}[V_L]_{(3-\chi(L))/2}\,=\,s(D)-1-\,\#\es\{\es(a,b)
\mbox{\ \ Seifert circles}\es:\,\mbox{
there is a crossing joining $a$ and $b$}\,\}\,.
\]
\end{theorem}

In other words, if $S_D$ is the reduced Seifert graph, then
$(-1)^{n(L)}[V_L]_{(3-\chi(L))/2}=b_1(S_D)$, $b_1$ being
the first Betti number.

\begin{corr}
For a positive diagram $D$, $b_1(S_D)$ is an invariant
of the link represented by $D$. \qed
\end{corr}

Note, that for the non-reduced Seifert graph $\bar S_D$
of $D$, $b_1(\bar S_D)=1-\chi(D)=1-\chi(L)$ is also a link invariant.

\begin{corr}
For a positive link $L$ of $n(L)$ components,
$(-1)^{n(L)}[V_L]_{(3-\chi(L))/2}\ge 0$, and this coefficient is
$0$ iff $L$ is fibered. \qed
\end{corr}

Of course, this fiberedness condition is not very useful,
when a positive diagram of $L$ is given, since to decide
then about fiberedness is trivial. However, applied in the
opposite direction, it can prove that $L$ is not positive.
This happens sometimes in a quite non-trivial way, as shows
the following example.

\begin{figure}[htb]
\[
\begin{array}{c}
\epsfsv{4.6cm}{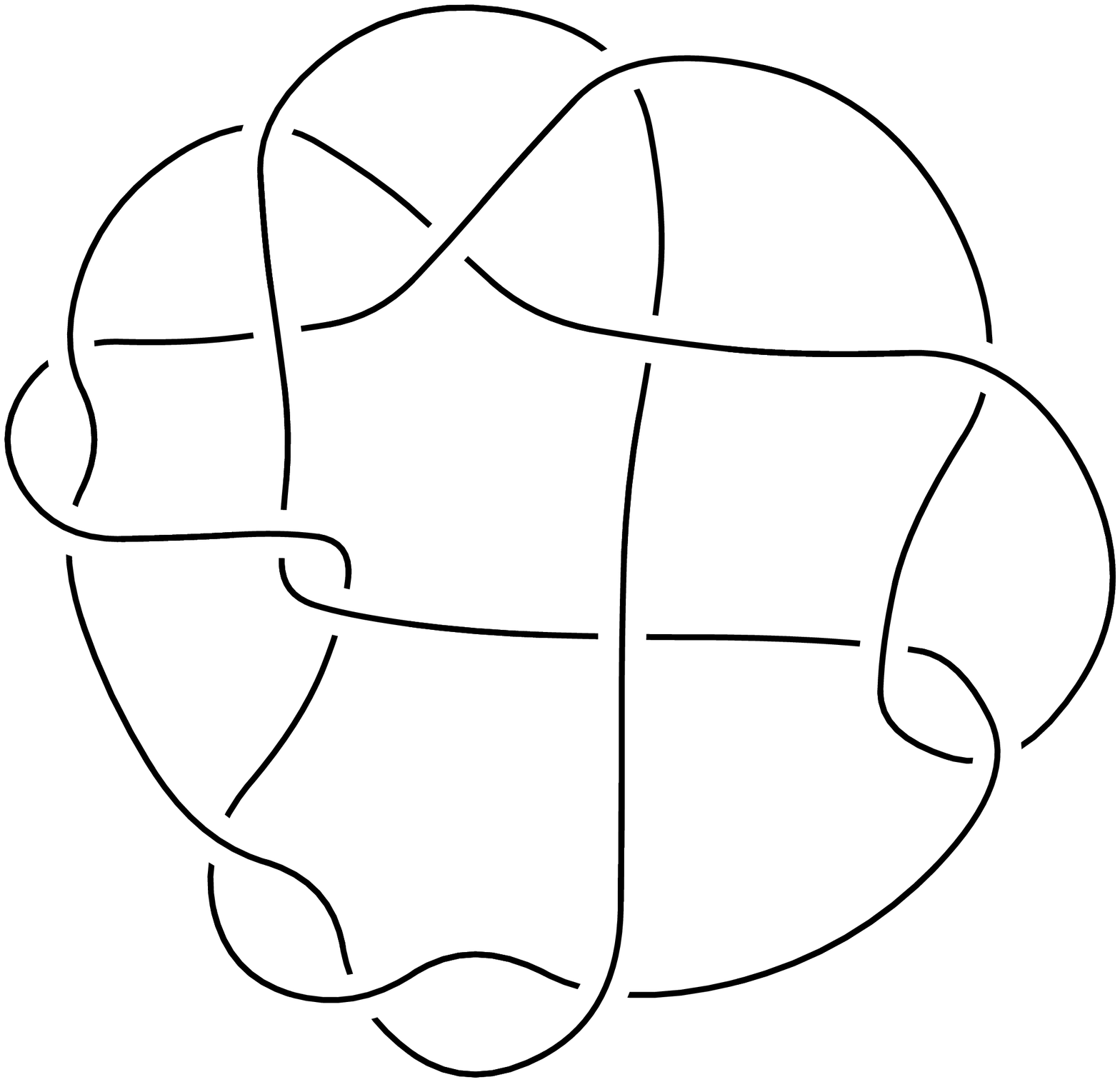} \\[21mm] 16_{1059787}
\end{array}
\]
\caption{\label{figP}}
\end{figure}

\begin{exam}
The knot $16_{1059787}$ in figure \reference{figP} satisfies
all conditions on positivity known about its $\nb$, $V$,
$P$ and $F$ polynomials. It seems useful to list all
properties that hold, even if they involve invariants we did
not consider here. See the given references for an accurate account.
(However, keep in mind that the conventions there differ from the ones
we use; for $F$ we conjugate in the $a$ variable.)
\begin{itemize}
\item $2\md V=\md_lP=\Md_mP=\Md\nb=\md_a F=4$ \cite{Cromwell,Fiedler,%
Yokota,Zulli},
\item $[P]_{m^4}(\sqrt{-l})$ and $\nb(z)$ are positive (that
is, all coefficients are non-negative) \cite{Cromwell},
\item $\tl P_i(l):=\sqrt{-1}^i[P]_{m^i}(\sqrt{-l})$ take only positive
values on $l\in(0,1)$ for $i=0,2,4$ \cite{MorCro}.
\item $[V]_{t^2}=1$ \cite{Fiedler,Zulli}.
\item $[F]_{a^4}(l)=[P]_{m^4}(l)$ \cite{Yokota}.
\item $[F]_{z^ka^k}\ge 0$ $\fa k\ge 0$ (that is, ``critical line''
polynomials are positive) \cite{Thistle2}.
\item $16_{1059787}$ does have diagrams (with canonical Seifert surface)
of genus $2$, so that $\tl g=g=\Md\Dl=2$ \cite{Cromwell}.
Here $\Dl$ is normalized so that $\Dl(t)=\Dl(1/t)$ and $\Dl(1)=1$.
\item The signature $\sg=4$ \cite{CG,apos}, so that by Murasugi's
inequality \cite{Murasugi2}, $g_s=g=2$ \cite{Rudolph,pos}.
\end{itemize}
However, now $[V]_{t^3}=0$, so that if $16_{1059787}$ is positive,
it must be fibered. But $[\Dl(t)]_{t^2}=2$ contradicts this property.
\end{exam}

\proof[of theorem \reference{cc}] This is proved as theorem
\reference{tht} using the bracket. The term $s(D)-1$
comes from the $A$-state, while for every pair of Seifert circles
joined by (at least) one crossing, a $-1$ comes from
an alternating sum of binomial coefficients coming from
states in which a $B$-splitting is applied at some
(non-empty) set of crossings linking $a$ and $b$. \qed

\begin{corr}\label{cr4}
Let $L$ be an almost positive link with an almost positive 
diagram $D$ such that there is no positive
crossing $q$ joining the same two Seifert circles as the
negative crossing $p$. Then
\[
\md V_L\,=\,\frac{1-\chi(D)}{2}\mbox{\ \ and\ \ }
\mcf V_L\,=\,(-1)^{n(L)-1}\,.
\]
\end{corr}

\proof Apply the skein relation at the negative crossing $p$
and use theorem \reference{cc} on $D_+$ and $D_0$ (they
have the same reduced Seifert graph). \qed

The following theorem is the key step needed to extend
Fiedler's result to almost positive links.

\begin{theorem}\label{thev}
Let $p$ be a crossing in a reduced special alternating diagram $D$
such that there is no crossing $q$ joining the same two Seifert circles
as $p$. Let $D_p$ be $D$ with $p$ smoothed out. Then
$\Dl_{D_p}(0)<\Dl_D(0)$, where $\Dl$ is the Alexander polynomial
normalized so that $\md\Dl=0$ and $\mcf\Dl=\Dl(0)>0$.
\end{theorem}

The proof will use the machinery of even valence graphs
\cite{MS}. We recall the basic notions from that paper.

\begin{defi}
The \em{join} (or \em{block sum}) `$*$' of two graphs is defined by
\[
\diag{6mm}{5}{2}{
  \picmultigraphics{2}{3 0}{\cycl{0 1}{1 0}{2 1}{1 2}}
  \picline{0 1}{2 1}
  \picputtext{2.5 1}{$*$}
}
\quad=\quad
\diag{6mm}{4}{2}{
  \picmultigraphics{2}{2 0}{\cycl{0 1}{1 0}{2 1}{1 2}}
  \picline{0 1}{2 1}
}
\]
This operation depends on the choice of a vertex in each one of
the graphs. We call this vertex the \em{join vertex}.

A \em{cut} vertex is a vertex, which disconnects the graph, when removed
together with all its incident edges. (A join vertex is always
a cut vertex.) Analogously a \em{$2$-cut} of $G$ is a pair of edges of
$G$ whose deletion disconnects $G$.
\end{defi}

\begin{defi}
A \em{cell} $C$ is the boundary of a connected component of the
complement of
a graph $G$ in the plane. It consists of a set of edges. If $p$ is
among these edges, then we say that $C$ \em{contains} $p$
or $p$ \em{bounds} $C$. By $G\sm C$ we mean the graph obtained from $G$
by deleting all edges in $C$.

A \em{cycle} $C$ is a graph $G$ is a set of edges $\{p_1,\dots,p_n\}$,
such that the pairs $(p_1,p_n)$ and $(p_i,p_{i+1})$ for $1\le i<n$
share a common vertex, and all these vertices are different. The plane
complement of a cycle in a planar graph has 2 components. The bounded
one we will call \em{interior} $\int(C)$ of $C$, and the unbounded
one \em{exterior} $\ext(C)$. (A cell is a cycle with one of
interior or exterior being empty, that is, containing no edges.)
\end{defi}

Before we make the next definition, first note, that the Seifert
graph $\bar S_D$ of any diagram $D$ is always planar(ly embeddable).
Namely $S_D$ is the join of the Seifert graphs corresponding to the
special diagrams in the Murasugi sum decomposition of $D$
along its separating Seifert circles, the join vertex
corresponding to the separating Seifert circle.
The join of planar graphs is planar, and 
if $D$ is a special diagram, then $\bar S_D$ has a natural
planar embedding (shrink the Seifert circles into vertices
and turn crossings into edges).

\begin{defi}
Assume for a special diagram $D$ that $\bar S_D$ is planarly embedded in
the natural way. Its dual is called the \em{even valence graph} $G_D$
of $D$ (as the name says, all its vertices have even valence).
Alternatively, $G_D$ is the checkerboard graph with vertices
corresponding to the non-Seifert circle regions of $D$.

A \em{canonical orientation} is an orientation of the
edges of $G_D$ so that all edges bounding a cell are oriented
the same way between clockwise and counterclockwise as seen from
inside this cell. (The canonical orientation is unique up to reversal
of orientation of all edges in a connected component of the graph.)
\end{defi}

\proof[of theorem \reference{thev}]
Consider the plane even valence graph $G_D$ associated to $D$. Then
$G_{D_p}=(G_D)_p$, where $G_p$ is $G$ with edge $p$ contracted.
Both $G$ and $G_p$ are connected by assumption. We shall assume from
now on that a canonical orientation is chosen in $G=G_D$, and
hence also on $G_p$.

By the matrix-tree theorem (see theorem 2 of \cite{MS}), we have
that $\Dl_D(0)=\mcf\Dl$ is the number of index-0 spanning rooted trees
of $G$, i.e. trees in which each edge, as oriented in $G$,
points towards the root of the tree. We will call such trees
also \em{arborescences}. Importantly, the number of arborescences
does not depend on the choice of root vertex. We will exploit this
property several times in the following.

Let $v_0$ be the  source, and $v_1$ the target of $p$ in $G$.
In $G_p$, $v_0$ and $v_1$ are identified to a vertex we call $v$.

By the proof of proposition 1, part 3), of \cite{MS}, we have that
\[
\#\,\{\mbox{\ index-0 sp. rooted trees with root $v$ in $G_p$\ }\}
\es=\es
\#\,\{\mbox{\ index-0 sp. rooted trees with root $v_1$ in $G$
containing $p$\ }
\}\,.
\]
Thus the statement of the theorem is equivalent to saying that
$G$ has an index-0 spanning rooted tree with root $v_1$ not
containing $p$. The assumption of the theorem in terms of
even valence graphs means that each edge of $G$ bounds a cell not
containing $p$, or equivalently, $p$ is in no $2$-cut of $G$.
In particular, both $v_0$ and $v_1$ have valence
at least $4$ in $G$.

It is easy to see that any plane even valence graph $G$ can
be built up from the empty one by adding directed cycles.
Moreover, if $G$ is connected, then we can achieve that
all intermediate graphs are connected (or more exactly
spoken, all their connected components except one are
trivial, i.e. an isolated vertex of valence 0).
Also, one can start the building-up with any particular
cycle in $G$.

Let $E$ be a cell (cycle with empty interior) in $G$ containing $p$.
We claim that then $\hG=G\sm E$ is still connected. There is a little
argument needed for this. We will explain that if a disconnected
graph $\hG$ is connected by adding a cell $E$, then each
edge in $E$ forms a $2$-cut with another edge in $E$. To see
this, first reduce the problem to $\hG$ having 2 components
$\hG_1$ and $\hG_2$. If $G$ has further components $\hG_3,\dots,
\hG_n$, one can connect them to $\hG_2$ by adding cells, and
a $2$-cut of edges in $E$ would still remain one if we undo this
connecting. It also is easy to see that one can assume there
are no valence-2-vertices of $E$ in $G$ (that is, each vertex of $E$
is attached to one of $\hG_1$ or $\hG_2$ in $G$). Then we show that
there are at most two edges of $E$ connecting $\hG_1$ and $\hG_2$.
Since $E$ is oriented, one can easily distinguish between
interior or exterior of $E$ depending on the (left or right) side
in orientation direction. If $\ge 4$ edges connect $\hG_1$ and $\hG_2$,
one must attach vertices of $\hG_1$ and $\hG_2$ from different
sides to $E$, and $E$ will not be a cell in $G$.

Let $E'$ be some
other cycle passing through $v_1$, such that $p\nin E'$.
(Such a cycle exists because $\val_G(v_1)>2$.) 

Then build up $G$ by adding cycles $E_n$, such that we start
with $E_1=E'$ and
finish with $E_z=E$, and all intermediate graphs $G_n$ are connected.
We construct successively in each $G_n$ an index-0
spanning rooted tree $T_n$ with root $v_1$, such that in the
final stage in $G_z=G$ the tree $T_z=T$ does not contain $p$.

In $G_1=E'$, fix the root to be $v_1$ and let $T_1$ consist
of all edges in $E'$ except the one outgoing from $v_1$.

\[
\diag{1cm}{4}{4}{
  \pictranslate{2 2}{
    \picPSgraphics{/xxi 0 D}
    \picfillgraycol{0}
    \picveclength{0.35}
    \picvecwidth{0.18}
    \picmultigraphics[rt]{8}{360 8 :}{
      \piclinewidth{/xxi xxi 1 + D xxi 2 ne {20}{30} ie}
      \picfilledcircle{2 0 polar}{0.08}{}
      \picvecline{2 0 polar}{2 360 8 : polar}{}
    }
    \picveclength{0.2}
    \picvecwidth{0.1}
    \picvecline{2.5 1.8}{2 360 8 : polar}
    \picputtext{2.1 1.4}{$p$}
    \picputtext{2 360 8 : polar 0.3 +}{$v_1$}
  }
}
\]

Now, given an index-0 spanning rooted tree $T_n$ of $G_n$,
we construct a index-0 spanning rooted tree $T_{n+1}$ of $G_{n+1}=
G_n\cup E_{n+1}$ as follows.

Let $w_1,\dots,w_k$ be the vertices of the cycle $E_{n+1}$
in cyclic order, so that $w_i$ and $w_{i+1}$ are connected
by a (directed) edge $p_i$. Then there is a non-empty set
$S\subset \{1,\dots,k\}$ such that for all $s\in S$, $w_s\in G_n$,
and $w_s$ is a trivial connected component (isolated vertex) in $G_n$
otherwise. Then add the following vertices to $T_n$ to obtain
$T_{n+1}$: for each $i,j\in S$ such that $(i,j)\cap S=\varnothing$ add
$\{\,p_m:\,m\in (i+1,j-1)\,\}$. Here $(i,j)$ is the interval of numbers
between (but not including) $i$ and $j$, meant w.r.t. the cyclic order
in $\bZ_k$.

Here is an example of a cycle $E_{n+1}$, in which the vertices in
$G_n$ are encircled, and the edges in $T_{n+1}\sm T_n$ thickened.

\[
\diag{1cm}{4}{4}{
  \pictranslate{2 2}{
    \picPSgraphics{/xxi 0 D}
    \picfillgraycol{0}
    \picveclength{0.35}
    \picvecwidth{0.18}
    \picmultigraphics[rt]{8}{360 8 :}{
      \piclinewidth{30}
      \picfilledcircle{2 0 polar}{0.08}{}
      \piclinewidth{/xxi xxi 1 + D xxi [3 5 6 8] x in {17}{30} ie}
      \piccircle{2 0 polar}{xxi [1 2 4 7] x in {0.26}{0.02} ie}{}
      \picvecline{2 0 polar}{2 360 8 : polar}{}
    }
  }
}
\]

Then $T_{n+1}$ is an index-0 spanning rooted tree with root $v_1$
in $G_{n+1}$.

It remains to see why $p\nin T_z=T$. For this note that $E=E_z\ni p$
is added last, and $\val_G(v_0),\val_G(v_1)\ge 4$, so that
$v_0,v_1\in G_{z-1}$. \qed

\begin{corr}\label{crt}
If $D$ is an almost positive diagram with negative crossing $p$ such
that there is no (positive) crossing $q$ joining the same two Seifert
circles as $p$, then $\md\Dl_D(t)=1-\chi(D)$, where $\Dl$ is
normalized so that $\Dl(t)=\Dl(1/t)$ and $\Dl(1)=1$.
In particular, the canonical Seifert surface associated to $D$
is of minimal genus.
\end{corr}

\proof Apply the skein relation for $\Dl$ at the negative crossing
to obtain the result for special diagrams. Then use the multiplicativity
of $[\Dl(D)]_{(1-\chi(D))/2}$ under Murasugi sum of diagrams
\cite{Murasugi2} to obtain the general case. \qed

This corollary improves the result of Hirasawa
\cite{Hirasawa} (theorem 2.1) stating that this Seifert
surface is incompressible.

Now we have all the preparations together to obtain the extension
of Fiedler's result.

\begin{theorem}\label{theo5}
If $L$ is an almost positive link, then $\md V_L=\ds\frac 12
(1-\chi(L))$, and $\mcf V_L=(-1)^{n(L)-1}$. 
\end{theorem}

\proof Let $D$ be an almost positive diagram of $D$ with negative
crossing $p$ and canonical Seifert surface $S$. One can easily reduce
the proof to the situation that $D$ is connected. Distinguish then two cases.
\def\labelenumi{\alph{enumi})}
\def\theenumi{\alph{enumi})}
\begin{enumerate}
\item 
There is a (positive) crossing $q$ joining the same two Seifert
circles as $p$.  By theorem \reference{tht} we must show that
\[
\frac{1-\chi(L)}{2}=\frac{1-\chi(D)}{2}-1\,.
\]
Clearly, $\ds\frac{1-\chi(L)}{2}\le \frac{1-\chi(D)}{2}$, and
by Bennequin's inequality $\ds\frac{1-\chi(L)}{2}\ge \frac{1-\chi(D)}{2}
-1$. 

Thus assume that $\ds\frac{1-\chi(L)}{2}=\frac{1-\chi(D)}{2}$, i.e.
$S$ is a minimal genus surface. By \cite{Gabai2}, then this
is true for the Murasugi summand of $S$, which is
the canonical Seifert surface associated to an almost positive
(or almost alternating) special diagram. However, by assumption
this surface is clearly not of minimal genus, a contradiction.

\item There is no such crossing $q$. Then we must show by corollary
\reference{cr4} that $\ds\frac{1-\chi(L)}{2}=\frac{1-\chi(D)}{2}$, i.e.
$S$ is a minimal genus surface. This follows again from
\cite{Gabai}, using corollary \reference{crt}. \qed
\end{enumerate}

\subsection{Skein polynomial and Morton's inequalities}

The results on the Jones polynomial and their proofs allow also
some applications to the skein polynomial \cite{HOMFLY}.

First, we can identify two more coefficients of the polynomial of
some positive links.

\begin{corr}
If $L$ is a 
fibered positive link of $n(L)$ components, then 
\[
[P_L]_{l^{1-\chi(L)}m^{-1-\chi(L)}}=(-1)^{n(L)}(1-\chi(L))\,,
\]
and if $L$ is prime and braid positive, then
\[
[P_L]_{l^{1-\chi(L)}m^{-3-\chi(L)}}=(-1)^{n(L)-1}
\,\chi(L)\,\frac{\chi(L)+1}{2}.
\]
\end{corr}

\proof Murasugi and Przytycki showed in \cite{MP} that $[P_D]_{m^{1-
\chi(D)}}$ is multiplicative under Murasugi sum. (That 
\begin{eqn}\label{ineqcg}
\Md_mP_D\le 1-\chi(D)
\end{eqn}
was shown by Morton \cite{Morton}.)
Since any positive diagram of a fibered positive link decomposes
as Murasugi sum of connected sums of $(2,\,.\,)$-torus links, we
have for any fibered positive link $L$,
\[
[P_L]_{m^{1-\chi(L)}}=l^{1-\chi(L)}\cdot (-1)^{n(L)-1}\,.
\]
Now apply corollary \reference{corr1} and the conversion 
\eqref{PtoV}\,.\qed

\begin{rem}
A formula for the first of the coefficients in the corollary
can be written for an arbitrary positive link using theorem
\reference{cc} instead of corollary \reference{corr1}.
\end{rem}

The proof of theorem \reference{theo5} can also be applied for $P$
and $\Dl$, normalized so that $\Dl(t)=\Dl(1/t)$ and $\Dl(1)=1$.

\begin{theorem}\label{theo5*}
If $L$ is an almost positive link, then 
\begin{eqn}\label{str}
2\Md\Dl(L)\,=\,\Md_mP(L)\,=\,1-\chi(L)\,.
\end{eqn}
\end{theorem}

\proof
That $2\Md\Dl_L\,=\,\Md_mP(L)$ follows from a well-known argument.
From \eqref{PtoDl} clearly $2\Md\Dl_L\le \Md_mP(L)$.
Assume $2\Md\Dl_L<\Md_mP(L)$.
This means that in the calculation of $\Dl_L$ by some (arbitrary)
skein resolution tree (see \cite[\S 1]{Cromwell}) contributions
of terminal nodes of the tree occur, which affect $\Md_lP_L$, but
cancel when substituting $l=i$ (as in \eqref{PtoDl}). However, in
the proof of \cite[corollary 2.2]{Cromwell}, Cromwell argued that
for an almost positive link diagram, one can choose the skein
resolution tree so that no contributions cancel. Thus we are left
to show the second equality in \eqref{str}.

Consider the two cases in the proof of theorem \reference{theo5}.

If $q$ shares its Seifert circles with another (positive) crossing in
$D$, then one of the special Murasugi sum components $D'$ of $D$,
can be reduced to a diagram $D''$ with $\chi(D'')>\chi(D)$. Thus
by \eqref{ineqcg}, 
\[
\Md_mP(D')\,=\,\Md_mP(D'')\le 1-\chi(D'')<1-\chi(D')\,,
\]
so that
$[P_{D'}]_{m^{1-\chi(D')}}=0$. Then by \cite{MP} the same holds for $D$.
Since we know that 
\begin{eqn}\label{str2.5}
1-\chi(L)=-1-\chi(D)
\end{eqn}
from the proof of theorem \reference{theo5}, the inequality 
\begin{eqn}\label{str3}
\Md_mP(L)\le 1-\chi(L)
\end{eqn}
follows.
Now use that, as a consequence of \cite[proposition 21]{LickMil},
for an arbitrary link $L$, $\md_lP\le \Md_mP$. From Morton's
inequalities \cite{Morton} we have then for an almost
positive diagram $D$ of $L$
\begin{eqn}\label{str4}
-1-\chi(D)\,=\,w(D)-s(D)+1\,\le\,\md_lP(D)\,\le\,\Md_mP(L)\,.
\end{eqn}
Now \eqref{str4} and \eqref{str2.5} show that the inequality
\eqref{str3} is exact.

If $q$ does not share its Seifert circles with another crossing in
$D$, then combining \cite{Morton} and the argument in the proof of
theorem \reference{theo5}, we have
\[
1-\chi(D)\,=\,2\Md\Dl(D)\,\le\,\Md_mP(D)\,\le\,1-\chi(D)\,=\,1-\chi(L)
\,,
\]
so we have the equality in \eqref{str}. \qed

\begin{rem}
Note that the equality $2\Md\Dl=1-\chi$ is a property
of purely algebraic-topological nature, and has been
examined for many links long before the skein polynomial
was around. Thus the usage of the skein polynomial in its
proof here is somewhat surprising. It seems essential, though,
and also explains why this equality cannot be proved so
for other classes. (Indeed, for example it is not always true
for almost alternating links.)
\end{rem}

As a consequence using theorem \reference{theo5}, we have

\begin{corr}
If $L$ is an almost positive link, then $\Md\Dl_L=\md V_L$. \qed
\end{corr}

This solves conjecture 5.2 of \cite{apos}, which stated this
property for almost positive knots. Using the work therein, we have

\begin{corr}
There are only finitely many almost positive knots with the
same Alexander polynomial. \qed
\end{corr}

Further we have

\begin{corr}
If $L$ is an almost positive link, then $\md_lP(L)\le 1-\chi(L)$.
\end{corr}

\proof Use again the mentioned consequence of
\cite[proposition 21]{LickMil}. \qed

This is another special case of Morton's conjectured inequality
\cite{Morton2} (disproved now in \cite{posex_bcr} for arbitrary links).
There is, though, much experimental evidence that we have
in fact equality in Morton's inequality.

\begin{question}\label{q2.1}
Is for any almost positive link $L$, $\md_lP(L)=1-\chi(L)$?
\end{question}

Note that in the one case in the proof of \eqref{str}, we obtained this
desired equality, namely when the almost positive diagram $D$ is
not of minimal genus. Latter property is to be
understood so that its canonical Seifert surface does not
realize the (Seifert) genus of $L$, i.e. $\chi(D)>\tl\chi(L)$.
Question \reference{q2.1} is thus related to the
question: Does any almost positive link $L$ have an almost positive
diagram $D$, which is not of minimal genus?

As we later found, the answer to this question is negative, and a
counterexample is the knot $!12_{1930}$ (which nevertheless satisfies
$\md_lP=1-\chi$). It is displayed in figure 8 of \cite{pos} (and occurs
also later in this paper as $L_4$ in the proof of corollary \ref{capo}).
Beside its obvious two almost positive diagrams (considered also
in the proof below), there are no other (reduced) ones. The proof of
this fact will be presented elsewhere, as it requires, apart from some
computation, several tools (developed in \cite{pos,apos,gen2}), that
go beyond the scope of the present paper.

The opposite situation to the last question is not less interesting,
in particular because positive diagrams are always genus-minimizing.

\begin{question}\label{q2}
Does any almost positive link $L$ have an almost positive diagram
$D$ of minimal genus? 
\end{question}

A positive answer to this question
will show that Morton's inequality for $\tl\chi(L)$,
\begin{eqn}\label{xxx}
\Md_mP(L)\,\le\,1-\tl\chi(L)
\end{eqn}
(which is a direct consequence of \eqref{ineqcg}) is sharp. It would
not be a surprise,
as knots with strict inequality are hard to find. So far two
methods apply: unity root values of $V$ \cite{gen2} and Gabai's
foliation algorithm \cite{Gabai3} to show that in fact
$\Md_mP(L)<1-\chi(L)$ \cite{posex_bcr}. Latter option seems
unlikely to work for almost positive links, and former option
requires considerable extension of the calculations. Out of the $\approx
4500$ non-alternating prime knots $K$ of $\le 16$ crossings with
$\Md_mP(K)\le 4$, in \cite{gen2} we obtained 28 such knots with
$4=\Md_mP(K)<2\tl g(K)$ using values of the Jones and $Q$ polynomial
at roots of unity (and one further undecided case). An easy check
shows that none of these 28 knots is almost positive.

\section{Almost positive diagrams with canonical fiber Seifert surfaces}

The even valence graphs can be used to give a description of
almost positive diagrams whose canonical Seifert surfaces are
fiber surfaces. The restriction to such surfaces is suggestive,
since in general the decision of fiber property of a link or a
surface may be difficult, even although both algebraic and geometric
methods are known. Our result is closely related to the result for
almost alternating diagrams due to Goda--Hirasawa--Yamamoto \cite{GHY}.
Our main motivation here was in fact to use the present (and
quite different) tools to extend and simplify the proof of their
criterion. We succeed almost completely, with the exception that we
cannot recover combinatorially the fact (see their proposition 5.1)
that instead of general Murasugi sum decomposability of the fiber into
Hopf bands in part \reference{it1} of the theorem we will state below
we have in fact stronger plumbing decomposability. On the other hand,
we show in part \reference{it3} that the fiberedness condition for the
Alexander polynomial is exact. Due to the copious ways to calculate
the Alexander polynomial, this makes the fiberedness property even
easier to detect than by the classification result \reference{it4}
for such diagrams by itself. (Our version of this result is also
more explicit than in the form given in \cite{GHY}.) In the next
section we will give examples showing (together with the examples
in \cite{GHY}) that one cannot extend the result much further.

In the following $\Dl$ will be normalized so that
$\Dl(t)=\Dl(t^{-1})$ and $\mcf\Dl>0$.

\begin{theorem}\label{th5}
Let $D$ be a connected almost positive link diagram with
canonical Seifert surface $S$. Then the following conditions
are equivalent:
\def\labelenumi{\arabic{enumi})}
\def\theenumi{\arabic{enumi})}
\begin{enumerate}
\item \label{it1} $S$ decomposes under iterated Murasugi sum (not
necessarily plumbing) completely into Hopf bands (of 1 full twist).
\item \label{it2} $S$ is a fiber surface.
\item \label{it3} $2\Md\Dl(D)=1-\chi(D)$ and $\mcf\Dl(D)=1$.
\item \label{it4} Decompose $D$ along its \em{separating} Seifert
circles (Seifert circles with non-empty interior and exterior) as
Murasugi sum of special diagrams, and those special diagrams
into prime factors. Then all these prime factors are
special alternating diagrams of $(2,n)$-torus links
(parallelly oriented), except for one, which after
reductions of the type
\begin{eqn}\label{y}
\diag{7mm}{2}{2}{
  \picPSgraphics{0 setlinecap}
  \pictranslate{1 1}{
    \picrotate{-90}{
      \lbraid{0 -0.5}{1 1}
      \lbraid{0 0.5}{1 1}
      \pictranslate{-0.5 0}{
      \picvecline{0.03 .95}{0 1}
      \picvecline{0.97 .95}{1 1}
    }
  }
  }
}
\quad\lra\quad
\diag{7mm}{1}{1}{
    \picmultivecline{0.18 1 -1.0 0}{0 0}{1 1}
    \picmultivecline{0.18 1 -1.0 0}{0 1}{1 0}
}
\es
\end{eqn}
becomes an almost positive special diagram of the following
forms:
\def\theenumi{\arabic{enumi}.}
\def\labelenumii{\alph{enumii})}
\def\theenumii{\alph{enumii})}
\begin{enumerate}
\item \label{itA} a special diagram whose (even valence)
checkerboard graph $G$ can be obtained as follows:

Take a chain of circles of positive edges
\begin{eqn}\label{3*}
\diag{1cm}{5}{1}{
  \ellipse{0.5 0.5}{0}{0.5 0.3}
  \ellipse{1.5 0.5}{0}{0.5 0.3}
  \ellipse{4.5 0.5}{0}{0.5 0.3}
  \picellipsearc{2.5 0.5}{0.5 0.3}{90 -90}
  \picellipsearc{3.5 0.5}{0.5 0.3}{-90 90}
  \picputtext{3 0.5}{$\dots$}
  \vrt{1 0.5}
  \vrt{2 0.5}
  \vrt{4 0.5}
}
\end{eqn}
and attach to it from outside a cell (cycle with empty exterior)
with one negative edge, which joins interior points of the
two 2 outermost loops in \eqref{3*}. (The negative edge corresponds
to the crossing to be switched.) E.g. see figure \reference{fig8}.

\begin{figure}
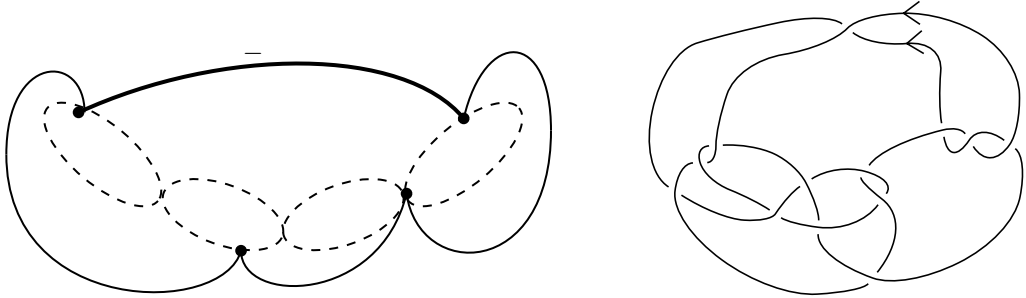

\[
\diag{8mm}{8}{6}{
  {\piclinedash{0.15 0.15}{0.25}
  \ellipse{1 3}{-40}{1.2 0.5}
  \ellipse{3 2}{-20}{1.05 0.5} 
  \ellipse{5 2}{20}{1.05 0.5} 
  \ellipse{7 3}{40}{1.2 0.5} 
  }
  {\piclinewidth{15}
   \piccurve{7 3.6}{6 4.8}{3 4.8}{0.6 3.7}
  }
  \piccurve{0.7 3.7}{0.7 4.7}{-0.6 4.7}{-0.6 3}
  \piccurveto{-0.6 0.3}{3 0.3}{3.3 1.4}
  \piccurveto{3.3 0.5}{5.6 0.5}{6.05 2.35}
  \piccurveto{6.25 0.9}{8.45 0.9}{8.45 3.4}
  \piccurveto{8.45 5.1}{7.35 5.1}{7 3.6}
  \vrt{7 3.6}
  \vrt{3.3 1.4}
  \vrt{6.05 2.35}
  \vrt{0.6 3.7}
  \picputtext{3.5 4.7}{$-$}
}
\rx{1.5cm}
\epsfsv{4.3cm}{t1-aal2}
\]
\caption{A checkerboard graph and its diagram illustrating
case \reference{itA} of theorem \reference{th5}.\label{fig8}}
\end{figure}

or (non-exclusively)

\item \label{itB}
a diagram of a $(2,2,\dots,2)$-pretzel link (at least two `$2$'s),
oriented to be special, with one crossing changed.
\end{enumerate}
\end{enumerate}
\end{theorem}

\proof In the following graph pictures, we assume graphs to
be canonically oriented, but do not draw edge orientation if
it is not necessary. The edge of the (only) crossing $p$ of negative
checkerboard sign will be distinguished by being drawn as a thickened
or dashed line.

$\ref{it4}\,\Ra\,\ref{it1}$. The reverse of the move on \eqref{y} preserves
the property the canonical Seifert surface to be a fiber, as it
corresponds to plumbing of a Hopf band. That the canonical Seifert
surfaces of the diagrams in \ref{itB} are fibers is easy to see.
(The diagram $D'$ obtained from $D$ by removing the trivial clasp
is a connected sum of Hopf links, its surface is clearly a fiber,
which is unique, and $\chi(D')=\chi(D)$.) For \ref{itA} we remark
that each of the graphs described turn into 
\[
\diag{1cm}{1}{0.5}{
  \picfillgraycol{0}
  \picfilledcircle{0 0.25}{0.08}{}
  \picfilledcircle{1 0.25}{0.08}{}
  \piccurve{0 0.25}{0.3 0.4}{0.7 0.4}{1 0.25}
  \piccurve{0 0.25}{0.3 0.1}{0.7 0.1}{1 0.25}
  \piccurve{0 0.25}{0.3 0.7}{0.7 0.7}{1 0.25}
  {\piclinewidth{20}
  \piccurve{0 0.25}{0.3 -0.2}{0.7 -0.2}{1 0.25}
  }
  \picputtext{0.5 -0.3}{$-$}
}
\]
under repeating the operation
\[
\diag{1cm}{4.5}{2}{
  \picveclength{0.4}\picvecwidth{0.14}
  \pictranslate{0.5 0}{
  \picfillgraycol{0}
  \picfilledcircle{1 1}{0.08}{}
  \picfilledcircle{3 1}{0.08}{}
  {\piclinewidth{14}
   \picline{1.8 -0.5 x}{1 1}
  }
  \picline{0.2 -0.5 x}{1 1}
  \picline{3 1}{4 1.3}
  \picline{3 1}{4 0.7}
  \piccurve{1 1}{1.6 1.7}{2.3 1.7}{3 1}
  \piccurve{1 1}{1.6 0.3}{2.3 0.3}{3 1}
  {\piclinedash{0.15 0.15}{0.25}
   \picellipsearc{1 1}{1.5 1}{25 -25}
  }}
}
\qquad\lra\qquad
\diag{1cm}{2}{2}{
  \picfillgraycol{0}
  \picfilledcircle{1 1}{0.08}{}
  {\piclinewidth{14}
   \picline{1.6 0 x}{1 1}
  }
  \picline{0.4 0 x}{1 1}
  \pictranslate{-2 0}{
    \picline{3 1}{4 1.3}
    \picline{3 1}{4 0.7}
  }
}
\]
(contracting a double edge), with the dashed line having two properties.
First, it is an arc passing through edges whose total sign sum is
0 (in our case the negative edge $p$ and one other, positive, edge), 
and second, it has a single vertex
and no complete edge in (at least) one of its interior or
exterior. This corresponds to the move on diagrams
\begin{eqn}\label{z}
\diag{1cm}{6}{3}{
  \picrotate{-90}{
    \lbraid{-0.75 1.25}{1.5 2.5}
    \lbraid{-0.75 3.75}{1.5 2.5}
    \piclinewidth{15}
    \picellipsearc{-0.75 5.2}{0.6 1.2}{220 320}
    \picellipsearc{-0.75 -0.2}{0.6 1.2}{40 140}
    \picellipse{-0.75 2.5}{0.6 1.0}{}
  }
  {\piclinedash{0.15 0.15}{0.25}
   \picellipsearc{2.5 1.55}{3.2 0.9}{-60 -120}
  } 
  \picputtext{2.5 2.7}{$\gm$}
}
\qquad\lra\qquad
\diag{1cm}{2}{1.5}{
  \picfillgraycol{0}
  \picrotate{-90}{
    \picellipsearc{-0.75 -0.6}{0.7 1.8}{40 140} 
    \picellipsearc{-0.75 3.6}{0.7 1.8}{220 320} 
    \piclinewidth{15}
    \picellipsearc{-0.75 3.2}{0.6 1.2}{220 320}
    \picellipsearc{-0.75 -0.2}{0.6 1.2}{40 140}
  }
}
\end{eqn}
on diagrams. The dashed line $\gm$ must pass through
interior of Seifert circles and crossings only, such that
the total writhe of these crossings is 0, and must have a
single non-Seifert circle region, and no complete
Seifert circle in (at least) one of its interior or exterior. 

Then this is a Hopf plumbing, the Hopf band
being obtained by thickening $\gm$ into a strip,
and taking the union with the two half-twisted
strips and one Seifert circle on the left of \eqref{z}.
(The first condition on $\gm$ is needed to ensure the correct
twisting, while the second one is needed to have the
Hopf band being separated by a sphere from the rest of the
surface after deplumbing.)

$\ref{it1}\,\Ra\,\ref{it2}\,\Ra\,\ref{it3}$ are well-known, so
it remains to show the real result $\ref{it3}\,\Ra\,\ref{it4}$.

As the minimal coefficient of the Alexander polynomial, when its degree
is equal to $(1-\chi)/2$, is multiplicative under
Murasugi sum \cite{Murasugi2}, we
need to consider only the almost alternating special Murasugi summand.

For this we consider the canonically oriented even valence graph $G=G_D$ of $D$ and
turn to the said in the proof of theorem \reference{thev}.
The condition $2\md\Dl=1-\chi$ implies that each edge of $G$ bounds
a cell not containing $p$, $p$ being the edge in $G=G_D$ of
negative checkerboard sign. In particular, both $v_0$ and
$v_1$ have valence $\ge 4$. 

Let $E_1$ and $E_2$ be the cells containing $p$. Then $E_1\cap E_2=
\{p\}$, since $p$ is in no $2$-cut of $G$, and $G\sm E_1$ and
$G\sm E_2$ are connected, by the argument in the
proof of theorem \reference{thev}. 

Now consider the condition $\mcf\Dl(D)=1$. It means that
there is only one index-0 spanning rooted tree with root $v_1$
not containing $p$. If $G\sm E_1$ has two different
index-0 spanning rooted trees with root $v_1$, then
by the construction in the proof of theorem \reference{thev},
we could extend them to index-0 spanning rooted trees
of $G$ with root $v_1$ not containing $p$, which would clearly
still be different. Thus $G\sm E_1$ has only one index-0 spanning
rooted tree (with root $v_1$ or any other fixed vertex). Then, by
part 5) of theorem 3 of \cite{MS}, $G\sm E_1$ is a join
of chains (topologically, a bouquet of circles).
\begin{eqn}\label{(7)}
\diag{1cm}{4}{3.5}{
  \drop{1.5 0.8}{-10}{0.8 1.2}
  \ddrop{1.5 0.8}{-160}{0.4 0.4}
  \drop{1.7 2}{30}{0.9 0.9}
  \ddrop{1.7 2.3}{80}{0.4 0.4}
  \ddrop{1.5 2.85}{-10}{0.4 0.7}
  \ddrop{1.1 2.6}{80}{0.5 0.7}
  \drop{1.5 0.8}{-65}{0.7 1.2}
  \ddrop{2.5 1.0}{-160}{0.4 0.7}
  \ddrop{2.5 1.0}{-95}{0.4 0.7}
  \drop{2.4 1.4}{-20}{0.6 0.9}
  \ddrop{2.5 2.2}{0}{0.5 0.7}
}
\end{eqn}
Assume w.l.o.g. that $G$ is embedded so that the exterior
$\ext(E_1)$ of $E_1$ is the unbounded component. 
Since adding $E_1$ must remove all cut vertices (our diagram
is prime by assumption), $E_1$ must touch interior points of
all circles $L_i$ with only one cut vertex. (Here `interior'
is meant with the circle considered 1-dimensional with
boundary being the cut vertex/vertices.) We call these $L_i$
\em{leaf} circles; in \eqref{(7)} they are drawn dashed.

Also, since the exterior of $E_1$ is the unbounded component,
cut vertices coming from attaching circles within others
\[
\diag{1.2cm}{1.3}{1.5}{
  \drop{0.7 0}{30}{1.2 1.2}
  \drop{0.7 0.7}{100}{0.4 0.4}
}
\]
cannot be removed by adding $E_1$, so assume there are no such
inner circles. Thus we have a picture like
\[
\diag{1cm}{4}{6}{
  \drop{1.5 0.8}{10}{0.7 1.0}
  \drop{1.5 0.8}{-50}{0.6 1.1}
  \drop{1.2 1.7}{20}{0.7 1.1}
  \drop{1.25 2.0}{-45}{0.6 0.9}
  \drop{2.0 2.6}{-35}{0.6 0.9}
  \drop{2.5 3.4}{0}{0.7 1.0}
  \drop{0.8 2.7}{20}{0.6 1.4}
  \drop{0.3 4.0}{0}{0.5 0.6}
  \drop{2.4 1.5}{-55}{0.6 0.9}
  \drop{3.2 2.0}{-40}{0.6 0.9}
  \drop{3.8 2.7}{-25}{0.5 0.8}
  \piclinewidth{25}
  \picputtext{1.0 5.4}{$p$}
  \picputtext{3.2 5.7}{$\ext(E_1)$}
  \picputtext{1.3 4.3}{$\int(E_2)$}
  \piccurve{4.0 3.4}{4.0 5.3}{2.5 5.3}{2.5 4.4}
  {\piclinedash{0.15 0.1}{0.05}
  \piccurveto{2.0 5.3}{1.0 5.6}{0.3 4.6}
  }
  \piccurveto{0.3 5.6}{-0.4 5.6}{-0.4 3.7}
  \piccurveto{-0.4 2.2}{0.3 1.8}{0.7 2.3}
  \piccurveto{0.5 1.9}{0.9 1.3}{1.1 1.4}
  \piccurveto{0.7 0.9}{0.7 0.3}{1.5 0.3}
  \piccurveto{2.6 0.3}{2.6 0.9}{2.4 1.5}
  \piccurveto{2.9 1.0}{3.4 1.2}{3.1 1.7}
  \piccurveto{3.9 1.3}{4.5 1.9}{4.1 2.9}
  \piccurveto{4.7 3.3}{4.5 3.8}{4.0 3.4}
}
\]

{}From now on, let us remove valence-2-vertices (we call this
operation \em{unbisection}) and consider
only the topological type of the tree
\begin{eqn}\label{**}
\diag{1cm}{2}{1}{
  \picline{0 0.5}{2 0.5}
  \vrt{1 0.5}
}
\qquad\lra\qquad
\diag{1cm}{1.5}{1}{
  \picline{0 0.5}{1.5 0.5}
}\quad.
\end{eqn}
This move on graphs corresponds to the reversed move
\eqref{y}. (Note that $\val_G(v_{0,1})\ge 4$ by assumption,
so that both edges on the left of \eqref{**} are positively
signed.)

The \em{way} between two leaf circles $L_1$ and $L_2$ is made up of
those circles bounding disks whose interior is passed by a path from an
interior point of the disk bounded by $L_1$ to an
interior point of the disk bounded by $L_2$. Hereby we require that
this path is passing only through
interior points of disks bounded by loops and cut vertices,
each such vertex being passed at most once.
\[
\diag{8mm}{9}{4}{
  \pictranslate{0 -1}{
  \ellipse{0.8 2.95}{-30}{1.4 0.6}
  \ellipse{3 2}{-20}{1.05 0.6} 
  \ellipse{5 2}{20}{1.05 0.6} 
  \ellipse{7 3}{40}{1.2 0.6} 
  \vrt{2 2.3}
  \vrt{4 1.7}
  \vrt{6 2.3}
  \vrt{0 3.5}
  \drop{2.8 2.7}{40}{0.4 0.9}
  \drop{2.8 2.7}{-20}{0.5 1.1}
  \drop{3.1 3.8}{0}{0.5 0.9}
  \drop{7 3.7}{-10}{0.7 1.1}
  \drop{7.5 2.6}{-120}{0.7 1.1}
  \drop{7.5 2.6}{-60}{0.5 1.4}
  \drop{8.7 3.3}{-20}{0.5 1.2}
  \vrt{9.03 4.2}
  {\piclinedash{0.15 0.1}{0.05}
  \piccurve{0 3.5}{0.6 3}{1.3 2.5}{2 2.3}
  \piccurveto{2.7 2.1}{3.3 1.7}{4 1.7}
  \piccurveto{4.5 1.7}{5.3 2}{6 2.3}
  \piccurveto{6.7 2.6}{7.3 2.5}{7.5 2.6}
  \piccurveto{7.7 2.7}{8.4 3.1}{8.7 3.3}
  \piccurveto{8.8 3.4}{8.9 3.8}{9.03 4.2}
  }
  \picputtext{0.5 4.3}{$L_1$}
  \picputtext{9.03 4.9}{$L_2$}
  }
}
\]

Now use that $G\sm E_2$ must also be a join of loops (or
bouquet of circles). We claim that either
\def\labelenumi{\alph{enumi})}
\def\theenumi{\alph{enumi})}
\begin{enumerate}
\item \label{it2A}
$p$ must touch interior points of two different leaf loops
$L_1$, $L_2$ and all other vertices of $E_1$ touch only 
interior points or cut vertices belonging to circles \em{on the
way} between $L_1$ and $L_2$ (as on the left of figure
\reference{fig8})

or
\item \label{it2B}
$p$ touches two points on the same loop $L$ (interior
points or the cut vertex), and $E_1$ touches points only of the
same loop.
\end{enumerate}

Assume neither \reference{it2A} nor \reference{it2B} holds.
We derive a contradiction showing that $G\sm E_2$ has at
least 2 arborescences. Observe that by unbisections \eqref{**}
and \em{separations}
\begin{eqn}\label{eq17}
\diag{1cm}{2}{2}{
  \picvecline{1 1}{0.5 2}
  \picvecline{0.5 0}{1 1}
  \picvecline{1 1}{0.2 0.5}
  \picvecline{0.2 1.5}{1 1}
  \picvecline{1.7 1.7}{1 1}
  \picvecline{1 1}{1.7 0.3}
  \vrt{1 1}
  \picputtext{1.5 0.8}{$e$}
}\quad\lra\quad
\diag{1cm}{2.7}{2}{
  \picvecline{1 1}{0.5 2}
  \picvecline{0.5 0}{1 1}
  \picvecline{1 1}{0.2 0.5}
  \picvecline{0.2 1.5}{1 1}
  \vrt{1 1}
  \picputtext{1.2 1.3}{$v'$}
  \pictranslate{0.7 0}{
    \picvecline{1.7 1.7}{1 1}
    \picvecline{1 1}{1.7 0.3}
    \picputtext{1 1.3}{$v$}
    \vrt{1 1}
    \picputtext{1.5 0.8}{$e$}
  }
}\,,
\end{eqn}
one can simplify $G$ to obtain
\begin{eqn}\label{7en}
\diag{1cm}{5}{5}{
  \pictranslate{2.5 2.5}{
    \picellipse{0 0}{1.5 0.75}{}
    {\piclinewidth{15}
     \picveclength{0.4}\picvecwidth{0.14}
     \picveccurve{1.5 210 polar 2 :}{-1 -2.7}{1 -2.7}{1.5 330 polar 2 :}
    }
    \picputtext{1.5 210 polar 2 : 0.3 + x 0.15 + x}{$v_0$}
    \picputtext{1.5 330 polar 2 : 0.3 + x 0.15 - x}{$v_1$}
    \vrt{1.5 210 polar 2 :}
    \vrt{1.5 330 polar 2 :}
    \vrt{0 0.75}
    \vrt{0 2.25}
    \picmultigraphics[S]{2}{-1 1}{
      \piccurve{0 0.75}{0.8 1.2}{0.8 1.8}{0 2.25}
      \opencurvepath{1.5 330 polar 2 :}{1.7 -0.8}{2.2 -0.6}{2.2 0.8}
        {2 1.5}{1.7 2.3}{0.9 2.5}{0 2.25}{}
    }
    \picputtext{0 -1.2}{$E_2$}
    \picputtext{2 -2}{$E_1$}
    \picputtext{0 -2.3}{$p$}
  }
}
\end{eqn}
(still with $E_1$ being the boundary of the infinite region),
in the way that $v_{0,1}$ are not involved in any of these moves.
Then removing $E_2$ and applying unbisections, one obtains a graph
$G_0$ with two vertices and an edge of multiplicity $4$, which has 2
arborescences.
It suffices now to show that \em{bisections} (reverses of unbisections)
and \em{deseparations} (reverses of separations \eqref{eq17}) do
not reduce the number of arborescences. For this we specify how to map
injectively arborescences of the original graph with root $v_1$
to arborescences of the resulting graph.

For a bisection creating a vertex $v\ne v_{0,1}$, add the outgoing
edge of $v$ to the arborescence, and so do with the incoming one, if
the original (bisected) edge was in the (original) arborescence.
The same argument, but without the restriction $v\ne v_{0,1}$, finds
2 arborescences (with root $v_1$ or any other vertex) of
\eqref{7en}, starting from those of $G_0$.

For a deseparation at least one of the two vertices $v,v'$ on
the right of \eqref{eq17} has valence $2$. Let $v$ be such a vertex.
The outgoing edge $e$ of $v$ is in any arborescence, since
$v\ne v_{0,1}$. Remove $e$ from the arborescence, and keep the
status of the other edges, while joining $v$ and $v'$.

Thus $G$ is one of the types \reference{it2A} and \reference{it2B}.

In case \reference{it2A} the assumption there are no
cut vertices in $G$ implies that $G\sm E_1$ is as in \eqref{3*},
and $p$ joins interior points on the two outermost circles.
Note, that the union of $L_{1,2}$ and all loops on the way
between them form a bouquet of type \eqref{3*}. Thus we arrive
at case \reference{itA} in theorem \reference{th5}.

In case \reference{it2B}, $G\sm E_1$ must be a single loop, and we
have a picture like this
\[
\diag{8mm}{3}{4.2}{
  \ellipse{1.5 2}{0}{1.0 1.6}
  \piclinewidth{19}
  {\piclinedash{0.15 0.1}{0.05}
  \piccurve{2.5 2.4}{3.9 3.1}{2.0 5.0}{1.5 3.6}
  }
  \piccurveto{1.0 5.4}{-1.1 2.5}{0.5 2}
  \piccurveto{-1.1 1.5}{0.6 -1.0}{1.1 0.5}
  \piccurveto{1.6 -0.4}{2.9 -0.1}{2.3 1.0}
  \piccurveto{3.3 1.3}{3.3 2.1}{2.5 2.4}
  \picputtext{2.6 3.5}{$p$}
}
\]
This is case \reference{itB} in theorem \reference{th5}.  \qed

What we have done allows to solve an open
problem in our previous work \cite{apos} on almost positive
knots and to prove

\begin{corr}\label{capo}
There exist almost positive knots of any genus $\ge 2$.
\end{corr}

It was shown in \cite{gen2}, that there are no almost positive knots of
genus $1$.

\proof Consider the $(3,3,\dots,3,-1)$-pretzel knots and links $L_n$.
(These are 2-component links if the number $n$ of `$3$'s is odd; in
this case we orient them so that the twists counted by the
`$3$'s are reverse.) As these diagrams of $L_n$ come from the
construction in part \reference{it4} of theorem \reference{th5},
$L_n$ are fibered (or see alternatively theorem 6.7 of \cite{Gabai4}).
Their diagrams reduce by one crossing to
$D_n=(-2\,-1,3,\dots,3)$ (one `$3$' less), which are almost positive
and of crossing number $3(1-\chi(D_n))=3(1-\chi(L_n))=3n$.
If $L_n$ were positive, by \cite[corollary 5.1]{Cromwell} they
would have crossing number $c(L_n)\le 2(1-\chi)$. 

To show that this is not the case, consider the crossing
number inequality of \cite{Kauffman,Murasugi,Thistle},
$c(L_n)\ge \spn V(L_n)$.

We know that $\ds\md V(L_n)=\frac{1-\chi}{2}$. On the other hand,
for $n=1-\chi>2$, $\Md V(L_n)$ is easy to determine, as the
diagram $D_n$ is $B$-semiadequate, and thus only the contribution
of the $B$-state $S_B$ ($\br{S_B\,:\,k}=k$ for any crossing $k$)
is relevant in \eqref{S}. By simple count of the loops one arrives
at $\Md V(L_n)=\myfrac{7}{2}(1-\chi)-2$, and thus
\[
c(L_n)\,\ge\,\spn V(L_n)\,=\,\Md V(L_n)-\md V(L_n)\,=\,3(1-\chi)-2\,
>\,2(1-\chi)\,.
\]
Thus $L_n$ is almost positive for $n\ge 3$. (For $n=1$ and $n=2$
one obtains the Hopf link and trefoil, resp.) \qed

\section{Some examples and problems}

\subsection{Showing almost positivity}

The problem to show that a certain link is almost positive,
but not positive, turned out very hard. All previously known
positivity criteria either are easily provable to extend to
almost positive links, or at least no examples are known
where they do not. Theorem \reference{theo5} is an addition
to that picture.

In \cite{apos} it was shown, in the case of knots,
that any almost positive knot has only finitely many
reduced almost positive diagrams. As the proof is constructive,
one can, in theory, decide (for knots) about almost positivity,
in the sense that for any knot one can write down a finite
set of almost positive diagrams, among which one would have to
check whether the knot occurs. However, this method is not
generally very efficient, except for a few knots of small genus.

Cromwell's estimate $c\le 2(1-\chi)$ for fibered homogeneous
links remains the only way known so far to circumvent these problems,
at least in certain cases. Using theorem \reference{th5}, we can
construct now a plenty of examples of almost positive fibered
link diagrams, which we can show to represent almost positive links
using that Cromwell's inequality is violated. 

However, this inequality will still not be violated in may cases,
and thus one may ask whether it can be improved.

Cromwell's estimate is trivially sharp for alternating (prime) links
(consider the rational links $222\dots2$) and composite positive links
(consider the connected sums of Hopf links). However, even for
prime positive links the inequality can not be improved much.

\begin{figure}[htb]
\[
\begin{array}{c@{\qquad}c}
\epsfsv{3.8cm}{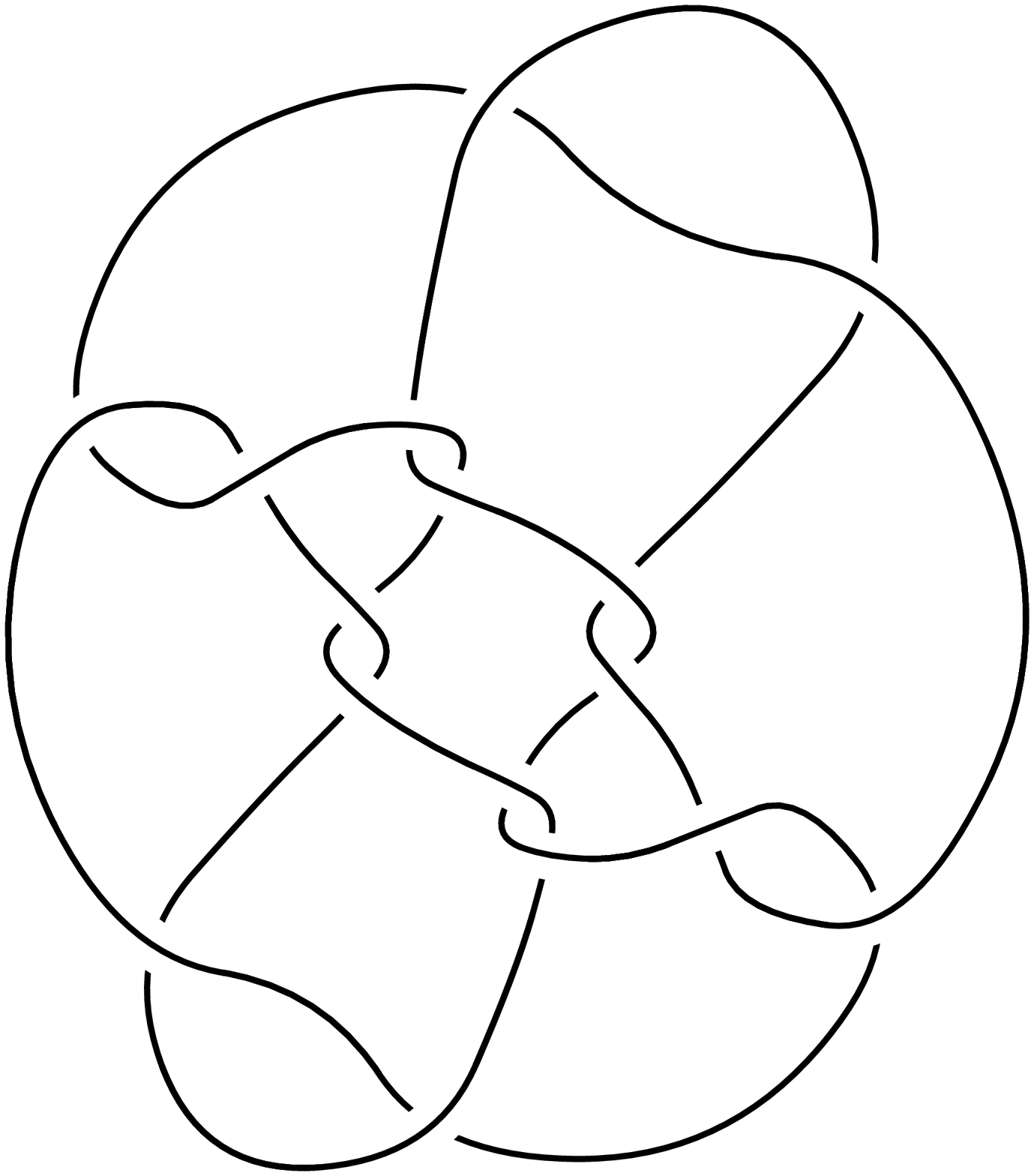} & \epsfsv{3.8cm}{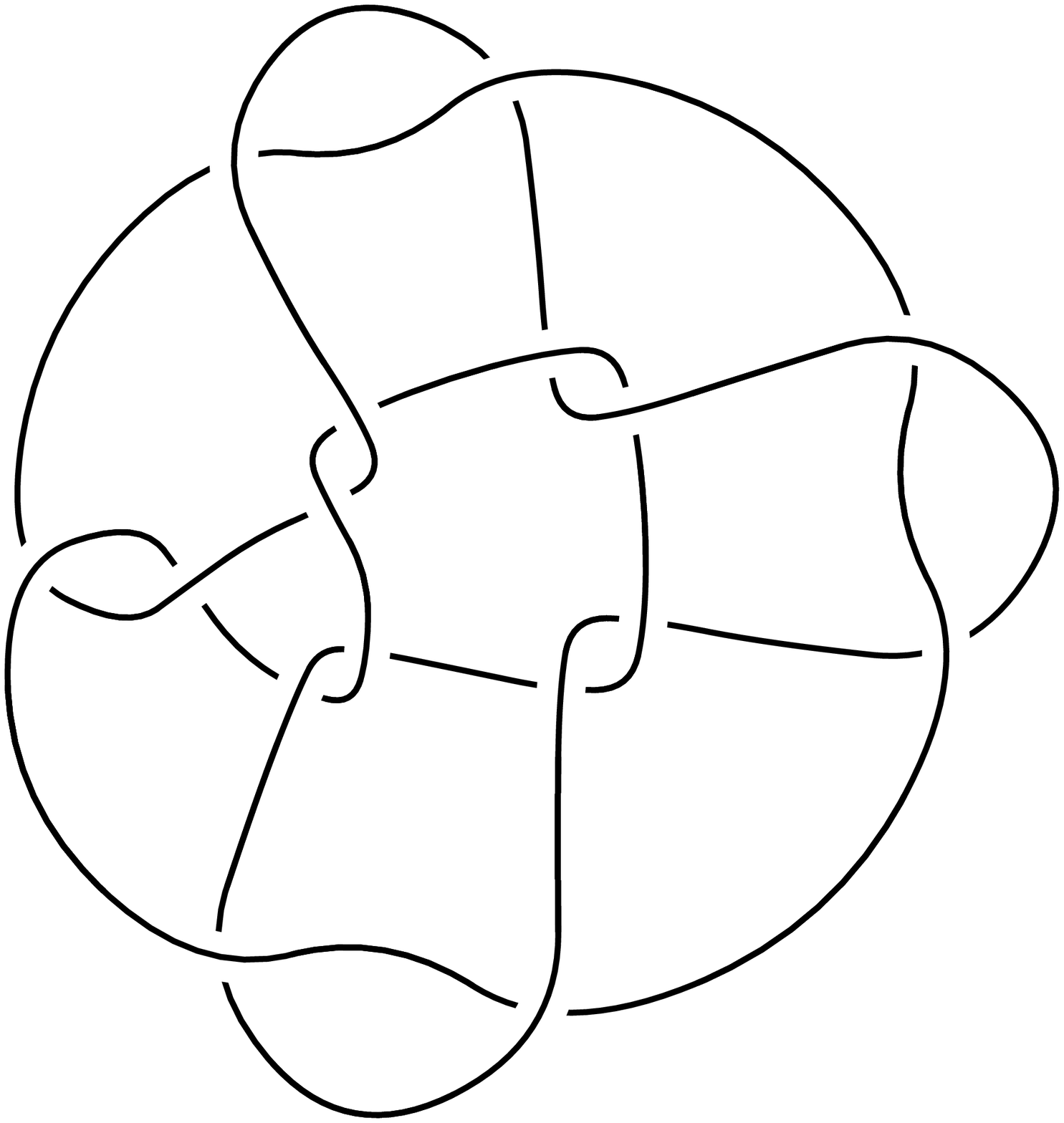} \\[21mm]
16_{1177344} & 16_{1243226}
\end{array}
\]
\caption{The two fibered positive knots of genus 4 and crossing number
16. (The diagrams here are chosen to be positive and reveal a plumbing
structure of the fiber surfaces. $16_{1243226}$ has also almost
positive 16 crossing diagrams, and $16_{1177344}$ even 2-almost
positive ones.)\label{fig4}}
\end{figure}

\begin{exam}
The $(2,2,\dots,2,-2,-2)$-pretzel link, oriented the
clasps counted by `$2$' to be reverse, and those counted by `$-2$'
to be parallel, has $c=2(1-\chi)-2$. The diagram is of minimal
crossing number as follows by considering linking numbers, and
Murasugi desums into connected sums of Hopf bands, thus the link is
fibered. (The link is also prime by \cite{KL}.)
\end{exam}

\begin{exam}
Even just considering
knots, there exist examples of genus 4 and crossing number 16,
$16_{1177344}$ and $16_{1243226}$. (For genus 3 the
maximal crossing number example is the knot $11_{550}$ of
\cite{posex_bcr} without a minimal positive diagram.)
Apparently these examples can be generalized to higher genera
(although the proof of minimality of the crossing number is
not straightforward). 
\end{exam}

Thus Cromwell's estimate seems rather sharp even in our case.
%

A different problem in this context is the position
of the class of almost positive links w.r.t. the hierarchy
\eqref{hie}.

\begin{question}
Is any almost positive link strongly quasipositive, or at least
quasipositive?
\end{question}

Some 2-almost positive links, like the figure-8-knot,
are not quasipositive. On the other hand, all almost positive
examples examined so far are strongly quasipositive.

\subsection{Detecting genus and fiberedness with
the Alexander polynomial}

{}From our results in the previous two sections, we have the following 

\begin{corr}\label{cx}
If a link $L$ has a connected almost positive (or almost alternating) diagram
(with canonical Seifert surface) of minimal genus, then
\def\labelenumi{\alph{enumi})}
\def\theenumi{\alph{enumi})}
\begin{enumerate}
\item $2\Md\Dl_L=1-\chi(L)$, and\label{cr1}
\item $L$ is fibered if and only if $\mcf\Dl_L= \pm 1$. \qed\label{cr2}
\end{enumerate}
\end{corr}

Unfortunately, we cannot decide about fiberedness if the
almost positive diagram is not of minimal genus. Many
almost positive knots seem to have almost positive 
minimal genus diagrams, but whether all have is unclear.
Coming back to the inequality $g(K)\le \tl g(K)$ in question \ref{q2},
it is known that almost alternating knots may fail to realize it
sharply. One of the two $\Dl=1$ knots of 11 crossings has genus two
\cite{Gabai3}, and is almost alternating by the verification in
\cite{Adams,Adams2}, while the calculation in \cite[example 11.1]
{LickMil} gives $\Md_mP=6$, so that by \eqref{xxx}, $\tl g=3$.
(A genus three canonical surface is not too hard to find.) This
knot thus does not have any diagram whatsoever of minimal genus.
Let us mention in contrast that among the 28 knots we found with
unsharp Morton inequality $4=\Md_mP(K)<2\tl g(K)$, none could be
identified as almost alternating (although there are not enough
tools to exclude it). However, there are several $\le 2$-almost
alternating knots, for example $15_{130745}$ (see figure 9 of
\cite{gen2}).

For almost positivity the problem to find knots with $\tl g>g$
seems much harder than for almost alternating.

\begin{question}
What is the minimal $k$ with a $k$-almost positive knot having
no diagram of minimal genus?
\end{question}

So far it seems likely that such knots with $k=4$ occur, but
even whether $k\le 3$ is unclear. Contrarily, there is a $2$-almost
positive knot, $16_{1337674}$, with unsharp inequality \eqref{xxx}.

Note that both statements in corollary \reference{cx} are true
for many (other) links, in particular all knots in Rolfsen's tables
\cite[appendix]{Rolfsen}. However, the following examples
show that the corollary does not extend much further.

\begin{figure}[htb]
\[
\begin{array}{c@{\qquad}c@{\qquad}c}
\epsfsv{3.8cm}{t1-2apos} & \epsfsv{3.8cm}{t1-2apos_13_6374} & 
\epsfsv{3.8cm}{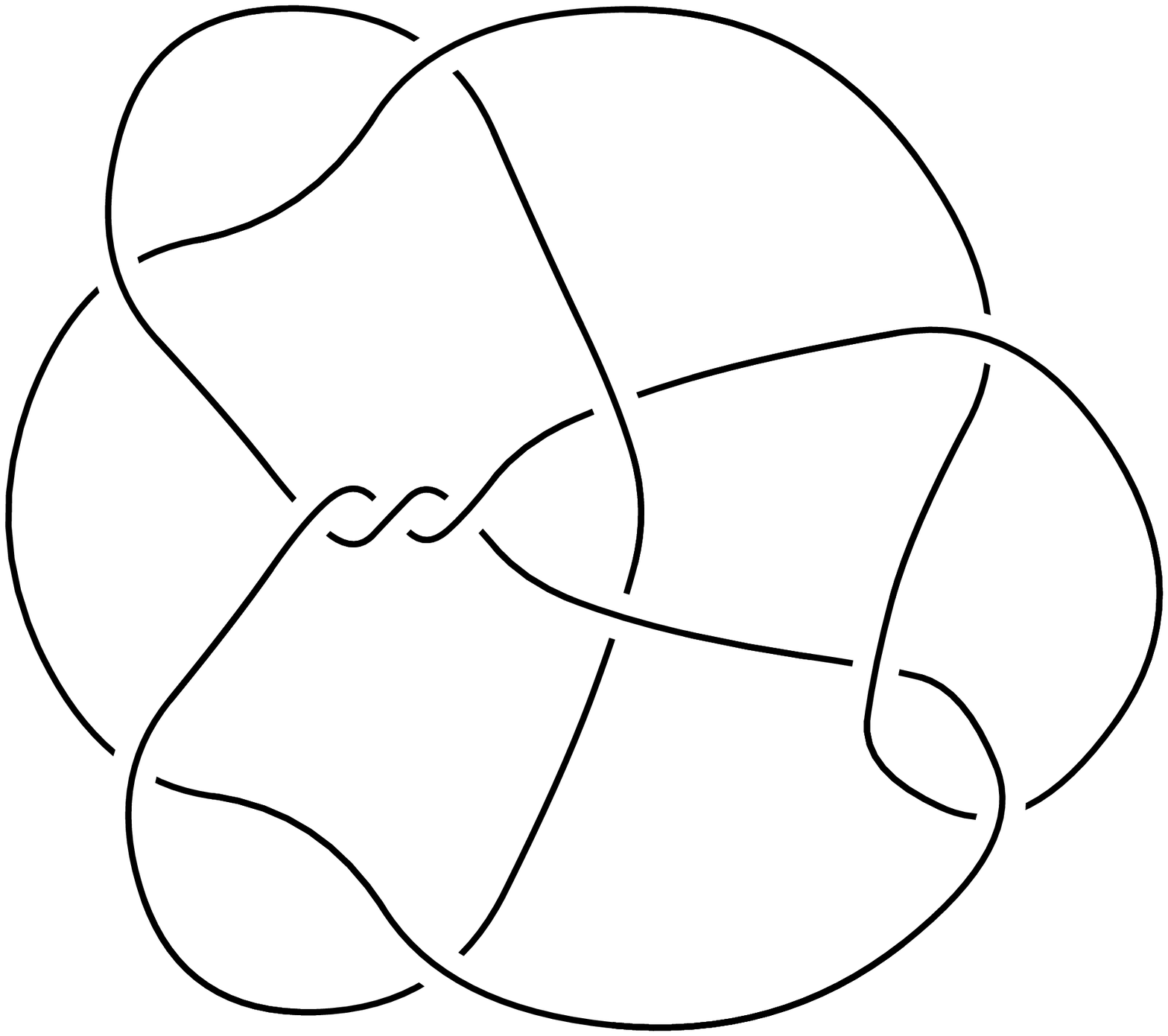}\\[21mm]
& 13_{6374} & 12_{1581}
\end{array}
\]
\caption{\label{f1}}
\end{figure}

\begin{exam}\label{ex1}
Consider the diagram on the middle of figure \reference{f1}.
It is another diagram of the previously encountered
knot $13_{6374}$ with unit Alexander polynomial.
It is 2-almost alternating, and its canonical Seifert surface
is of minimal genus (two), as can be shown by \cite{Gabai3}.
Thus neither of both criteria hold for 2-almost alternating
diagrams.
\end{exam}

\begin{exam}
The diagram on the right of figure \reference{f1}
depicts the knot $12_{1581}$ with Alexander polynomial
$\Dl=(\,2\ [-3]\ 2\,)$. It is a (special) 2-almost positive diagram
whose canonical Seifert surface is of minimal genus (again two).
Thus criterion \reference{cr1} in corollary \reference{cx}
is not true for 2-almost positive diagrams.
\end{exam}

So far I have no example of a 2-almost positive \em{knot} diagram for
criterion \reference{cr2}, but one can easily obtain a link
diagram.

\begin{exam}
Consider the diagram of $13_{6374}$ in example \reference{ex1}.
It has a single separating Seifert circle, whose interior
contains two crossings. By removing this interior (deplumbing
a Hopf band), one arrives at the link diagram on the left
of figure \reference{f1}. Its canonical Seifert surface is still
of minimal genus by \cite{Gabai2}, so that $1-\chi=3$, but one
calculates that $\Dl=t^{1/2}-t^{-1/2}$. 
\end{exam}

In all the above examples we showed a
surface not to be a fiber using that the Alexander polynomial
has too small degree. There are also examples, where the degree 
is maximal, and thus all conditions in corollary \reference{cx}
taken together still do not suffice to determine a fiber.

\begin{figure}[htb]
\[
\begin{array}{c@{\qquad}c}
\epsfsv{3.8cm}{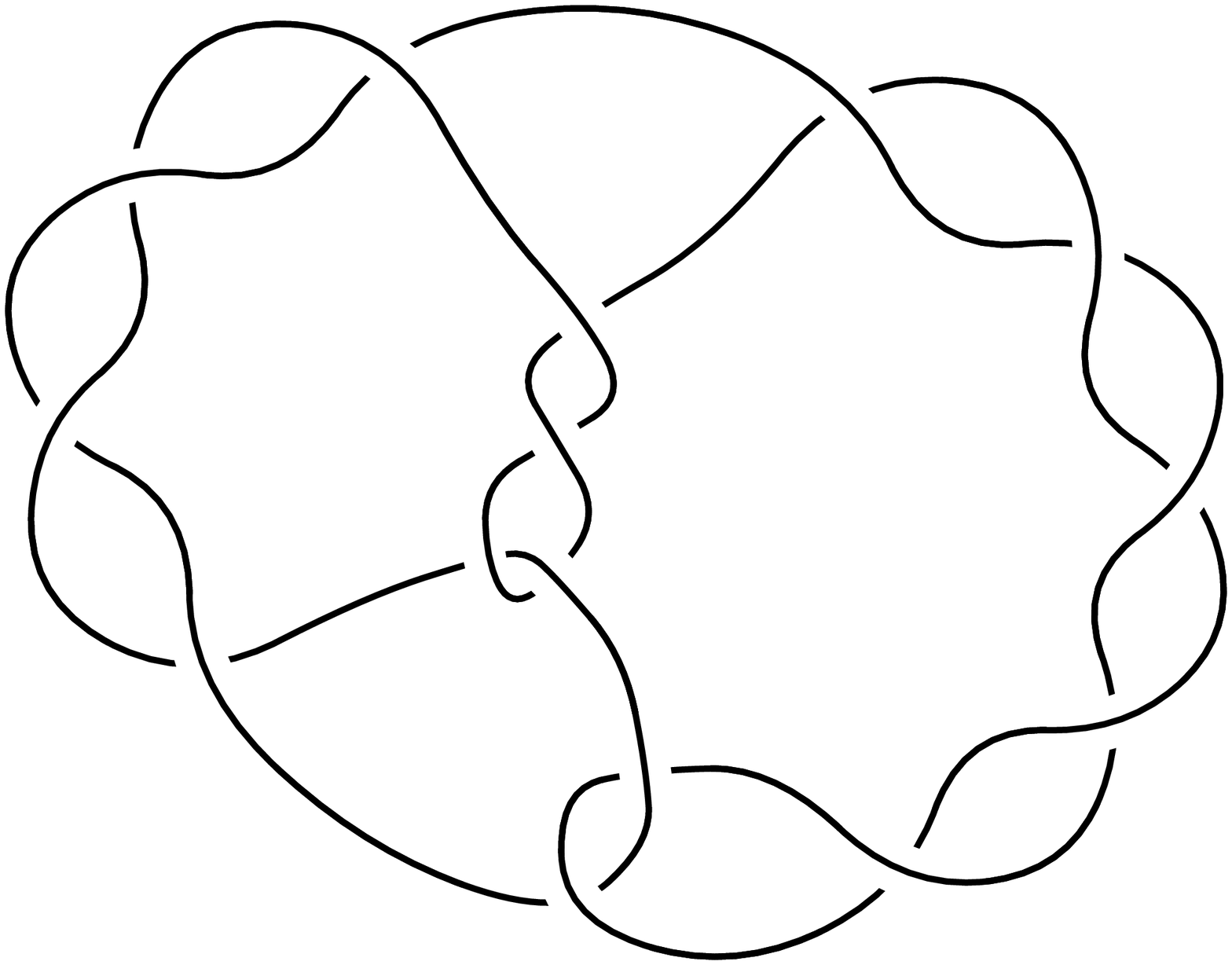} & \epsfsv{3.8cm}{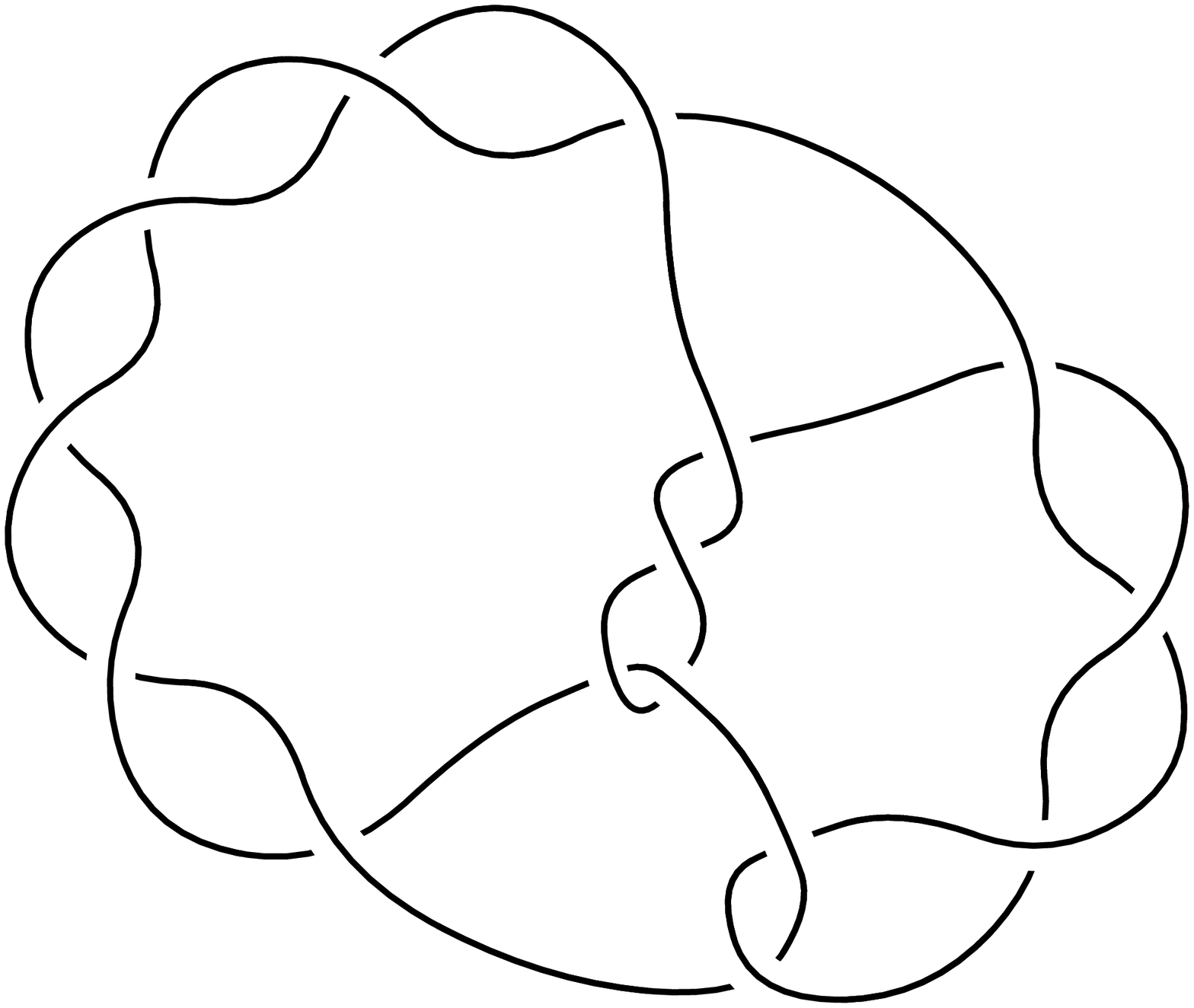} 
\\[21mm]
15_{120617} & 15_{159580}
\end{array}
\]
\caption{\label{f2}}
\end{figure}

\begin{exam}
The $(-2,4,6)$-pretzel link diagram, oriented to be special
(all clasps reverse), has $\Md\Dl=1-\chi=2$ and $\mcf\Dl=1$.
That its canonical surface is not a fiber follows from
\cite[theorem 6.7]{Gabai4} (Case 1). Using properly signed Hopf
plumbings, one can generate from this one many more examples of
2-almost alternating and/or 2-almost positive diagrams, in particular
(diagrams of) several genus two knots. Two such knots (for
2-almost positive diagrams) are the mirror images of $15_{120617}$
and $15_{159580}$, displayed in figure \reference{f2}. (These two knots
have been in fact found first, by a check in the tables, and the
pretzel link was obtained from them by deplumbings.)
\end{exam}

\begin{rem}
For genus one canonical surfaces of knots one needs
3-almost positive (and 3-almost alternating) diagrams to
have such a situation, the simplest example being the $(-3,5,5)$-pretzel
knot \cite{CT}. (It has the Alexander polynomial of the figure-8-knot.)
For 4-almost positive (or 4-almost alternating) diagrams,
even worse, one can use \cite[theorem 6.7]{Gabai4} to construct
diagrams differing by mutation, with the canonical surface of
the one being a fiber, and of the other one not.
\end{rem}

\noindent{\bf Acknowledgement.} I would wish to thank to M.\ Khovanov
for helpful remarks and discussions, in particular for suggesting to
me Theorem \reference{th1} and, implicitly, question \reference{qu1}.
The idea to distinguish in the study
of almost positive diagrams between having another crossing
joining the same pair of Seifert circles as the negative
one or not appears previously in Hirasawa's
paper \cite{Hirasawa}, and was observed even before that
in unpublished work of K.\ Taniyama.

{\small
\let\old@bibitem\bibitem

}
\end{document}